\documentclass{article}
\usepackage[centering,includeheadfoot,margin=2.5 cm]{geometry}

\usepackage{graphics,epsfig,enumerate}
\usepackage{amsfonts,amsmath,amssymb}
\usepackage{euscript,color,mathrsfs,accents}

\newtheorem{convention}{Convention}
\newtheorem{remark}{Remark}
\newtheorem{defin}{Definition}
\newtheorem{theorem}{Theorem}

\newtheorem{lemma}{Lemma}

\def\begproof{\noindent{\bf Proof: }}
\def\endproof{\par\rightline{\vrule height5pt width5pt depth0pt}\medskip}

\newcommand{\ubar}[1]{\underaccent{\bar}{#1}}    

\def\({\begin{eqnarray}}
\def\){\end{eqnarray}}
\def\[{\begin{eqnarray*}}
\def\]{\end{eqnarray*}}

\def\part#1#2{\frac{\partial #1}{\partial #2}}
 
\def\tot#1#2{\frac{\d #1}{\d #2}} 

\def\Norm#1{\left\| #1 \right\|}

\def\half{\frac{1}{2}}

\def\d{\,\mathrm{d}}
\def\R{\,\mathbb{R}}
\def\N{\mathbb{N}}
\def\E{\mathbb{E}}
\def\A{\mathbb{A}}
\def\e{\mathbf{e}}

\def\eps{\varepsilon}
\def\x{\mathbf{x}}
\def\vecv{\mathbf{v}}
\def\w{\mathbf{w}}
\def\u{\mathbf{u}}
\def\bxi{\boldsymbol{\xi}}

\def\eqref#1{$(\ref{#1})$}

\begin{document}

\title{On Cucker-Smale model with noise and delay\thanks{The research 
leading to these results has received funding from the European 
Research Council under the {\it European Community}'s Seventh Framework 
Programme ({\it FP7/2007-2013})/ ERC {\it grant agreement} No. 239870.}} 

\author{Radek Erban$^1$, 
\hskip 3mm 
Jan Ha\v{s}kovec$^2$  
\hskip 2mm 
and 
\hskip 2mm 
Yongzheng Sun$^3$
}

\maketitle

\footnotetext[1]{Mathematical Institute, University of Oxford, Radcliffe Observatory Quarter,
Woodstock Road, Oxford, OX2 6GG, United Kingdom; e-mail: erban@maths.ox.ac.uk. 
Radek Erban would like to thank the Royal Society for a University  Research Fellowship  
and the Leverhulme Trust for a Philip Leverhulme Prize.}
\footnotetext[2]{Mathematical and Computer Sciences and Engineering Division, 
King Abdullah University of Science and Technology, Thuwal 23955-6900, 
Kingdom of Saudi Arabia; e-mail: jan.haskovec@kaust.edu.sa.}
\footnotetext[3]{School of Sciences, China University of Mining and Technology, Xuzhou 221116,  
PR China; e-mail: yzsung@gmail.com. Yongzheng Sun acknowledges financial support from 
the China Scholarship Council (Grant No. 201308320087) and the National Natural Science 
Foundation of China (Grant No. 61403393).}

\begin{abstract}
A generalization of the Cucker-Smale model for collective animal 
behaviour is investigated. The model is formulated as a system of 
delayed stochastic differential equations. It incorporates two 
additional processes which are present in animal decision making, 
but are often neglected in modelling: (i) stochasticity (imperfections) 
of individual behaviour; and (ii) delayed responses of individuals 
to signals in their environment. Sufficient conditions for flocking 
for the generalized Cucker-Smale model are derived by using 
a suitable Lyapunov functional. As a byproduct, a new result 
regarding the asymptotic behaviour of delayed geometric Brownian 
motion is obtained. In the second part of the paper results of 
systematic numerical simulations are presented. They not only 
illustrate the analytical results, but hint at a somehow surprising 
behaviour of the system - namely, that an introduction of intermediate 
time delay may facilitate flocking.
\end{abstract}

\textbf{Keywords:} Cucker-Smale system, flocking, asymptotic behaviour, 
noise, delay, geometric Brownian motion.

\pagestyle{myheadings}
\thispagestyle{plain}
\markboth{R. Erban, J. Haskovec and Y. Sun}{On Cucker-Smale model with noise and delay}

\section{Introduction}\label{sec:Intro}
Collective coordinated motion of autonomous self-propelled
agents with self-organization into robust patterns
appears in many applications ranging from animal herding to the emergence 
of common languages in primitive societies~\cite{Sumpter}.
Apart from its biological and evolutionary relevance, collective phenomena
play a prominent role in many other scientific disciplines, such
as robotics, control theory, economics and social 
sciences~\cite{Carrillo-review, Vicsek-survey, Pareschi-Toscani-survey}.
In this paper we study the interplay of noise and delay on collective
behaviour. We investigate a modification of the well known Cucker-Smale 
model~\cite{CS1, CS2} with multiplicative noise and reaction delays.

We consider $N\in\N$ autonomous agents described by their phase-space 
coordinates $(x_i(t), v_i(t))\in\R^{2d}$, $i=1,2,\dots,N$,
where $x_i \equiv x_i(t) \in \R^{d}$ (resp. $v_i \equiv v_i(t) \in \R^{d}$) 
are time-dependent position (resp. velocity) vectors of the $i$-th agent. 
The governing equations are given as the following system of delayed
It\^{o} stochastic differential equations
\begin{eqnarray}
\d x_i &=& v_i \d t, 
\label{CS01} 
\\
\d v_i &=& \frac{\lambda}{N} \sum_{j=1}^N 
\psi_{ij}
\, 
\big( 
{\widetilde v_j}
-
{\widetilde v_i}
\big) \d t 
+ 
\frac{\sigma_i}{N} \sum_{j=1}^N  
\psi_{ij} 
\, 
\big( 
{\widetilde v_j}
-
{\widetilde v_i}
\big)
\d B_{i}^t,  
\label{CS02} 
\end{eqnarray}
where the delayed velocity is given by
$
\widetilde v_i (t) = v_i(t-\tau)
$
and $\tau\geq 0$ is the reaction delay. The parameters $\lambda>0$ 
and $\sigma_i \in \R$, $i=1,2,\dots,N,$ measure the alignment and noise 
strength, respectively, and $\d B_{i}^t$, $i=1,2,\dots,N$, are 
independent $d$-dimensional white noise vectors.
In general, the communication rates $\psi_{ij}$ 
are functions of the mutual distances $|x_i-x_j|$, however, in most of our paper we will
consider them as given functions of time satisfying certain assumptions.
The standard Cucker-Smale model~\cite{CS1, CS2} is a special case of 
equations~(\ref{CS01})--(\ref{CS02}) for $\sigma_i = 0$ and $\tau = 0$.
Our aim is to investigate equations (\ref{CS01})--(\ref{CS02}) for general
values of reaction delay $\tau$ and noise strength parameters
$\sigma_i$, $i=1,2,\dots,N.$ 

The Cucker-Smale model was introduced and studied in the seminal 
papers~\cite{CS1, CS2}, originally as a model for language evolution. 
Later the interpretation as a model for flocking in animals 
(birds) prevailed. In general, 
the term \emph{flocking} refers to the phenomena where autonomous agents reach 
a consensus based on limited environmental information and simple rules. 
The Cucker-Smale model is a simple relaxation-type model that reveals 
a phase transition depending on the intensity of communication between agents. 
Using $\sigma_i = 0$ and $\tau = 0$ in (\ref{CS01})--(\ref{CS02}), we can 
write the standard Cucker-Smale model as the following system of ordinary 
differential equations 
\begin{eqnarray}
\dot x_i 
&=& 
v_i \,,  
\label{CSm1}
\\
\dot v_i 
&=& 
\frac{\lambda}{N} \sum_{j=1}^N 
\psi_{ij} (v_j-v_i) \,,\qquad\qquad\mbox{for } i=1,2,\dots,N,
\label{CSm2}
\end{eqnarray}
where the dots denote the time derivatives. We note that the scaling 
by $N^{-1}$ in~(\ref{CSm2}) is significant to obtain a Vlasov-type 
kinetic equation in 
the mean-field limit $N\to\infty$, see, for example \cite{Tadmor-Ha}. The 
communication rates $\psi_{ij}$ introduced in~\cite{CS1, CS2} and most of 
the subsequent papers are of the form
\begin{equation} 
\label{commRate}
\psi_{ij}=\psi(|x_i-x_j|) \qquad\qquad\mbox{with } \;   
\psi(s) = \frac{1}{(1+s^2)^\beta},
\quad \mbox{where} \; \beta\geq 0.
\end{equation}
If $\beta < 1/2$, then the model exhibits the so-called 
\emph{unconditional flocking}, where for every initial configuration
the velocities $v_i(t)$ converge to the common \emph{consensus value}
$N^{-1} \sum_{i=1}^N v_i(0)$ as $t\to\infty$.
On the other hand, with $\beta\geq 1/2$ the flocking is \emph{conditional},
i.e., the asymptotic behaviour of the system depends on the value of $\lambda$
and on the initial configuration. 
This result was first proved in \cite{CS1, CS2} using tools from graph theory
(spectral properties of graph Laplacian), and slightly later reproved 
in~\cite{Tadmor-Ha} by means of elementary calculus. Another proof has been 
provided in~\cite{Ha-Liu}, based on bounding \eqref{CSm1}--\eqref{CSm2}
by a system of dissipative differential inequalities, and, finally, the 
proof of~\cite{CFRT} is based on bounding the maximal velocity.

Various modifications of the generic model \eqref{CSm1}--\eqref{CSm2} have 
been considered. For instance, the case of singular communication rates 
$\psi(s) = 1 / s^\beta$ was studied in~\cite{Ha-Liu, Peszek}.
Motsch and Tadmor~\cite{Motsch-Tadmor} scaled the communication rate 
between the agents in terms of their relative distance, so that their 
model does not involve any explicit dependence on the number of agents.
The dependence of the communication rate on the topological rather 
than metric distance between agents was introduced in \cite{Haskovec}.
The influence of additive noise in individual velocity measurements 
was studied in~\cite{Ha-Lee-Levy} and~\cite{TLY}.
Stochastic flocking dynamics with multiplicative white noises was 
considered in~\cite{Ahn-Ha}. Delays in information processing were 
considered in~\cite{Liu-Wu}, however, their analysis only applies to 
the Motsch-Tadmor variant of the model.

In this paper, we are interested in studying the combined influence of 
noise and delays on the asymptotic behaviour of the Cucker-Smale system. 
In particular, we derive a sufficient condition in terms of noise 
intensities $\sigma_i$ and delay length $\tau$ that guarantees flocking. 
Our analysis is based on a construction of a Lyapunov functional and an 
estimate of its decay rate. To prove our main results, we make an additional 
structural assumption about the matrix of communication rates which, 
loosely speaking, means that the communication between agents 
is strong enough.

The paper is organized as follows: In Section \ref{sec:Model}
we show how our model \eqref{CS01}--\eqref{CS02} is derived from 
the Cucker-Smale model \eqref{CSm1}--\eqref{CSm2} and define
what is meant by \emph{flocking}. Moreover, we consider a simplified version
of the model to provide an intuitive understanding of what qualitative
properties may be expected.
In Section~\ref{sec:Flocking} we derive a sufficient condition for flocking
in terms of the parameters $\lambda$, $\sigma_i$ and $\tau$,
based on a micro-macro decomposition and construction of a Lyapunov 
functional. Moreover, as a byproduct of our analysis, we provide 
a new result about the asymptotic behaviour of delayed geometric 
Brownian motion. Section~\ref{sec:Numerics} is devoted to a systematic 
numerical study of the model. First, we focus on simulation of delayed 
geometric Brownian motion, in particular, we study the dependence of
its asymptotic behaviour on the delay and noise levels. Then, we perform 
the same study for system~(\ref{CS01})--(\ref{CS02}). This leads to 
the interesting observation that, for weak coupling and small noise 
levels, an introduction of intermediate delays may facilitate flocking.
A systematic study of this effect concludes the paper.

\section{The stochastic Cucker-Smale model with delay}
\label{sec:Model}
In order to make the generic model~(\ref{CSm1})--(\ref{CSm2}) more 
realistic, we amend it with two additional features. First, we note 
that measurements in the real world are subject to errors and 
imprecisions that are typically modeled in terms of white noise.
In particular, we assume that the state (velocity) of agent $j$ 
measured by agent $i$ is given by the expression
\begin{equation} 
\label{omega_i}
\omega_{i;j} = v_j + \kappa_i (v_j - v_i) \d B_{i}^t,
\qquad \mbox{for} \quad i,j=1,2,\dots,N,
\end{equation}
where $\kappa_i\geq 0$ represents the imprecision of $i$'s measurement 
device, and $B_{i}^t$ are independent
identically distributed $d$-dimensional Brownian motions with zero mean
and the covariance relations
\begin{equation*} 
\E \left[ \! B_{i}^t \cdot B_{j}^s \right] = d \, \delta_{ij}\delta_{ts} t,
\qquad \mbox{for} \quad i,j=1,2,\dots,N \mbox{ and } t,s\geq 0,
\end{equation*}
with $\delta_{ij}=1$ iff $i=j$ and $\delta_{ij}=0$ otherwise,
and similarly for $\delta_{ts}$.

Note that the multiplicative structure of the noise term
ensures that $\omega_{i;i} = v_i$. Substituting $\omega_{i;j}$ given by 
(\ref{omega_i}) for $v_j$ in (\ref{CSm2}) and defining 
$\sigma_i:=\lambda\kappa_i$, we obtain the following system of 
stochastic differential equations (SDEs) for velocities
\begin{equation*} 
\d v_i 
= 
\frac{\lambda}{N} \sum_{j=1}^N \psi_{ij} (v_j-v_i) \d t
+ 
\frac{\sigma_i}{N} \sum_{j=1}^N  \psi_{ij} (v_j-v_i) \d B_{i}^t,
\qquad \mbox{for} \quad i=1,2,\dots,N,
\end{equation*}
with $\lambda$ a positive coupling strength.

The second amendment is the introduction of delays, motivated by the fact
that agents react to information received from their surroundings with some time lag.
However, we assume that information propagates instantaneously, so the delay
does not depend on the physical distance between agents.
For simplicity, we assume the reaction lag to be the same for all agents,
so that at time $t$ they react to information perceived at time $t-\tau$
for a fixed $\tau>0$.

\medskip

\begin{convention}
Throughout the paper, we denote by $v_i$ the quantity $v_i$ evaluated at time $t$,
i.e., $v_i=v_i(t)$, and by $\widetilde v_i$ the same quantity evaluated 
at time $t-\tau$, i.e., $\widetilde v_i = v_i(t-\tau)$.
We will also write $\x=(x_1,x_2,\dots,x_N)\in\R^{d\times N}$ 
$($resp. $\vecv=(v_1,v_2,\dots,v_N)\in\R^{d\times N}$$)$
for the vectors of locations $($resp. velocities$)$ of the agents.
\end{convention}

\medskip

In general, the communication rates $\psi_{ij}>0$ may be functions
of the mutual distances $|x_i-x_j|$. However, our analysis
is based on a certain structural assumption about the communication matrix
$(\psi_{ij})_{i,j=1}^N$ and the particular form of the dependence
on the mutual distances is irrelevant. Therefore,
we consider the rates $\psi_{ij}=\psi_{ij}(t)$ as given adapted
stochastic processes, so that \eqref{CS02} decouples from \eqref{CS01}.
Moreover, we assume that $\psi_{ij}$ are uniformly bounded,
\begin{equation} 
\label{PsiBoundSymm}
0 < \psi_{ij}(t) = \psi_{ji}(t) \leq 1,
\qquad\mbox{for } \; t\geq 0,
\quad i,\,j=1,2,\dots,N,
\quad\mbox{almost surely.}
\end{equation}
Thus, we finally arrive at the stochastic system of delayed differential
equations that we will study,
\( \label{CS0}
\d v_i &=& \frac{\lambda}{N} \sum_{j=1}^N 
\widetilde\psi_{ij}
\, 
\big( 
{\widetilde v_j}
-
{\widetilde v_i}
\big) \d t 
+ 
\frac{\sigma_i}{N} \sum_{j=1}^N  
\widetilde\psi_{ij} 
\, 
\big( 
{\widetilde v_j}
-
{\widetilde v_i}
\big)
\d B_{i}^t,  
\)
which is supplemented with the deterministic constant initial datum
$\vecv^0 \in \R^{d\times N}$,
\begin{equation}  
\label{IC}
v_i(t) \equiv v_i^0, \qquad\mbox{ for } \; t\in(-\tau,0],\quad i=1,2,\dots,N.
\end{equation}
Let us note that we interpret the noise term in $(\ref{CS0})$
in terms of the It\^{o} calculus \cite{Oksendal, Mao}.

\begin{theorem}\label{thm:ex}
The stochastic delay differential system~$(\ref{CS0})$ with initial
condition~$(\ref{IC})$ admits a unique global solution 
$\vecv=\vecv(t)$ on $[-\tau,\infty)$ 
which is an adapted process with 
$\E\left[\int_{-\tau}^T |\vecv(t)|^2 \d t\right] < \infty$ for all $T<\infty$,
i.e., a martingale.
\end{theorem}

\begproof
The proof follows directly from Theorem 3.1 and the subsequent remark on p. 157 of \cite{Mao}.
Indeed, (\ref{CS0}) is of the form
\[
    \d\vecv(t) = F(t,\vecv(t-\tau))\d t + G(t,\vecv(t-\tau))\d B^t
\]
for suitable functions $F$ and $G$.
In particular, the right-hand side is independent of the present state $\vecv(t)$,
so that the solution can be constructed by the method of steps. 
The second order moment is bounded on $(-\tau,T)$
because of the linear growth of the right-hand side of (\ref{CS0}) in $\widetilde\vecv$.
\endproof
\medskip

We now define the property of \emph{asymptotic flocking} for the solutions
of~(\ref{CS0})--(\ref{IC}).
\medskip

\begin{defin}\label{def:flocking}
We say that the system $(\ref{CS0})$ exhibits asymptotic flocking if
the solution $(\vecv(t))_{t\geq 0}$ for any initial condition $(\ref{IC})$
satisfies
\begin{equation*}
\lim_{t\to\infty} 
\Big| 
\E[v_i(t)] - \E[v_j(t)] 
\Big| = 0 \qquad\qquad\mbox{{\rm for all} } \; i,\, j =1,2,\dots,N,
\end{equation*}
where $\E[\cdot]$ denotes the expected value of a stochastic process.
\end{defin}

\subsection{Simplified case with $\psi\equiv 1$}
To get an intuition of what qualitative properties
we may expect from the solutions of \eqref{CS02},
we consider the case when the communication rate
is constant, i.e., $\psi_{ij}\equiv 1$;
in other words, we set $\beta=0$ in~(\ref{commRate}).
We also assume that $\sigma_i$ is equal to the same
constant $\sigma \in \R$ for all $i=1,2,\dots,N$, i.e. 
$\sigma_i \equiv \sigma$, and, moreover, that $v_i^0=v^0$
for some $v^0\in\R^d$ and all $i=1,2,\dots,N$.
Then, defining $V_c(t):=\frac{1}{N}\sum_{i=1}^N v_i(t)$,
we obtain
\[
\d V_c = 
\frac{\sigma}{N} \sum_{i=1}^N \left( {\widetilde V_c} 
- {\widetilde v_i} \right) \d B_{i}^t.
\]
Since, by assumption, $V_c(t) - v_i(t)\equiv 0$ for $t\in(-\tau,0]$,
we have $V_c(t)\equiv v^0$ for all $t\geq 0$.
Consequently, \eqref{CS0} decouples
into $N$ copies of the delayed SDE
\begin{equation}
\label{simplified}
\d w = - \lambda \, \widetilde w \, \d t - \sigma \, \widetilde w \, \d B^t,
\end{equation}
where we denoted $w:=v_i-v^0$ for any $i=1,2,\dots,N$.
We are not aware of any results concerning
the asymptotic behaviour of~\eqref{simplified}.
The method developed in~\cite{AMR} suggests that
$$
\lim_{t\to\infty} 
\E \left[|w(t)|^2 \right] = 0
\qquad
\qquad 
\mbox{if and only if}
\qquad
\qquad 
\int_0^\infty r_\lambda(t)^2 \d t < \frac{1}{\sigma^2},
$$
where $r_\lambda$ is the fundamental solution of the delayed ODE
\begin{equation}  
\label{ODE}
\dot w = -\lambda \, \widetilde w,
\end{equation}  
i.e., formally, $r_\lambda$ solves~(\ref{ODE}) subject to the 
initial condition $w(t) = \chi_{\{0\}}(t)$ for $t\in(-\tau,0]$.
The fundamental solution $r_\lambda$ can be constructed by the 
method of steps~\cite{Smith}, however, evaluation of its 
$L^2(0,\infty)$-norm is an open problem.
From this point of view, the analysis carried out in 
Section \ref{sec:Flocking} provides new and valuable information 
about the asymptotics of~(\ref{simplified}), see 
Section \ref{subsec:delayedGBM}. Let us note that setting 
$\tau=0$ in the above criterion recovers the well-known
result about geometric Brownian motion~\cite{Oksendal}:
the mean square fluctuation $\E[|w(t)|^2]$ tends to zero
if and only if $\sigma^2 < 2\lambda$.

Finally, for the convenience of the reader, we give an overview
of the qualitative behaviour of solutions to~(\ref{ODE}) 
with $\lambda>0$, subject to a constant nonzero initial 
datum (see, e.g., Chapter 2 of~\cite{Smith}):
\begin{itemize}
\item If $\lambda\tau\leq 1/e$, the solution monotonically 
converges to zero as $t\to\infty$, hence no oscillations occur.
\item If $1/e < \lambda\tau < \pi/2$, oscillations appear, however, 
with asymptotically vanishing amplitude.
\item If $\lambda\tau = \pi/2$, periodic solutions exist.
\item If $\lambda\tau > \pi/2$, the amplitude of the oscillations 
diverges as $t\to\infty$.
\end{itemize}
Hence, we conclude that the (over)simplified model \eqref{ODE},
corresponding to the delayed Cucker-Smale system with $\psi\equiv 1$ 
and no noise,
exhibits flocking if and only if $\lambda\tau < \pi/2$.
In the next Section we derive a sufficient condition for flocking
for the general model \eqref{CS0}.

\section{Sufficient condition for flocking}
\label{sec:Flocking}
In this section we derive a sufficient condition for flocking
in \eqref{CS0} according to Definition \ref{def:flocking}.
Our analysis will be based
on a construction of a Lyapunov functional
that will imply decay of velocity fluctuations for suitable parameter values.
However, we will have to adopt an additional structural assumption
on the matrix of communication rates $(\psi_{ij})_{i,j=1}^N$.

Before we proceed, let us shortly point out the mathematical difficulties
that arise due to the introduction of delay and noise 
into the Cucker-Smale system.
The ``traditional'' proofs of flocking of model~(\ref{CSm1})--(\ref{CSm2}), 
for instance~\cite{CS1, CS2, Tadmor-Ha, Ha-Liu}, rely on the monotone decay 
of the kinetic energy (velocity fluctations) of the form
$$
\tot{}{t} \sum_{i=1}^N |v_i|^2 
= 
- \frac{\lambda}{N} \sum_{i=1}^N\sum_{j=1}^N \psi_{ij} |v_i-v_j|^2 \leq 0.
$$
However, this approach fails if processing delays are introduced,
since for \eqref{CS02} without noise (i.e., all $\sigma_i=0$), we have
$$
\tot{}{t} \sum_{i=1}^N |v_i|^2 
= - \frac{\lambda}{N} \sum_{i=1}^N\sum_{j=1}^N 
\widetilde\psi_{ij} (v_i-v_j)\cdot(\widetilde v_i-\widetilde v_j).
$$
One then expects the product $(v_i-v_j)\cdot(\widetilde v_i-\widetilde v_j)$
to be nonnegative for $\tau>0$ small enough,
however, it is not clear how to prove this hypothesis.

The introduction of noise leads to additional difficulties
- in particular, the classical bootstrapping argument \cite{CS1, CS2, Ha-Liu}
for fluctuations in velocity fails in this case.
Similarly as in~\cite{Ha-Lee-Levy}, we circumvent this 
problem by adopting, in addition to the boundedness (\ref{PsiBoundSymm}),
a structural assumption about the matrix of communication rates.
We define the Laplacian matrix $\A(t)\in\R^{N\times N}$ by
\begin{equation} 
\label{Psi}
\A_{ij} 
:= -\psi_{ij}, \quad \mbox{ for } \; i\neq j,
\qquad \qquad
\A_{ii} := \sum_{j\neq i} \psi_{ij},
\end{equation}
and note that $\A$ is symmetric, diagonally dominant with non-negative 
diagonal entries, thus it is positive semidefinite and has real nonnegative 
eigenvalues. Due to its Laplacian structure,
its smallest eigenvalue is zero~\cite{CS1}.
Let us denote its second smallest eigenvalue (the Fiedler number) $\mu_2(t)$.
Our structural assumption is that there exists an $\ell > 0$ such that
\begin{equation} 
\label{fiedlervect}
\mu_2(t) \geq \ell > 0, 
\qquad \qquad \mbox{for } \; t > 0,\quad\mbox{almost surely.}
\end{equation}
This can be guaranteed for instance by assuming that the communication rates
are uniformly bounded away from zero, $\psi_{ij}(t) \geq \ubar{\psi} > 0$,
since there exists a constant $c>0$ such that $\mu_2(t) \geq c\ubar{\psi}$,
see Proposition 2 in~\cite{CS1}.

Moreover, we assume that the matrix of communication rates
is uniformly Lipschitz continuous in the Frobenius norm, in particular,
there exists a constant $L>0$ such that
\( \label{Lipschitz}
   \Norm{\A(t)-\A(t-\tau)}_F \leq L\tau \qquad\mbox{for } t\geq 0,
\)
where $\Norm{\cdot}_F$ denotes the Frobenius matrix norm.

To ease the notation and without loss of generality,
we will consider the one-dimensional setting $d=1$, i.e., 
$v_i(t) \in \R$ and $\vecv(t)\in\R^N$, where $N$ is the number of agents.
Then, with the definition~(\ref{Psi}), we put~(\ref{CS0}) 
into the form
\begin{equation} 
\label{CS3}
\d v_i 
= 
-\frac{\lambda}{N} (\widetilde\A\widetilde\vecv)_i \d t 
+ 
\frac{\sigma_i}{N} (\widetilde\A\widetilde\vecv)_i \d B^t_i,
\qquad \quad i=1,2,\dots,N.
\end{equation}
Our main result is the following.

\medskip

\begin{theorem}
\label{thm:flocking}
Let $\A$ be given by $(\ref{Psi})$ satisfying $(\ref{PsiBoundSymm})$, $(\ref{fiedlervect})$
and $(\ref{Lipschitz})$. Let the parameters $\lambda>0$ and
$\sigma^2_\mathrm{max}:=\max \{\sigma^2_1,\sigma^2_2,\dots,\sigma^2_N\}$ 
satisfy
\begin{equation}  
\label{flockingCond_0}
\sigma^2_\mathrm{max} < \lambda,
\end{equation}
then there exists a critical delay $\tau_c = \tau_c(\lambda,\sigma_\mathrm{max},L,\ell)>0$,
independent of $N$, such that for every $0 \leq \tau < \tau_c$
the system~$(\ref{CS3})$ exhibits flocking in the sense of Definition 
$\ref{def:flocking}$.

Moreover, if the matrix of communication rates $\A$ is constant,
i.e. $(\ref{Lipschitz})$ holds with $L=0$,
then $\tau_c$ is of the form
\begin{equation}  
\label{flockingCond}
\tau_c = \frac{1}{\lambda^2}
\left( - \, \sigma_\mathrm{max}^2 
+ 
\sqrt{\, \sigma_\mathrm{max}^4 
+ \frac{1}{12} \, (\lambda-\sigma_\mathrm{max}^2)^2} \right).
\end{equation}
\end{theorem}

\medskip

\begin{remark}\label{rem:noell}
The system $(\ref{CS3})$ with
constant communication matrix $\A$ can be seen
as a linearization of the system $(\ref{CS01})$--$(\ref{CS02})$
about the equilibrium $v_i\equiv v_0$ for $i=1,2,\dots,N$
with some $v_0\in\R$.
Note that in this case the formula~$(\ref{flockingCond})$ for the critical
delay $\tau_c$ does not depend on the 
particular value of $\ell$ in~$(\ref{fiedlervect})$.
\end{remark}

\subsection{Micro-macro decomposition}\label{subsec:micro-macro}
We introduce a micro-macro decomposition \cite{Tadmor-Ha, Ha-Lee-Levy}
which splits \eqref{CS3} into two parts: macroscopic, that describes
the coarse-scale dynamics, and microscopic, that describes the
fine-scale dynamics. The macroscopic part for the
solution is set to be 
the mean velocity $V_c(t)$,
\begin{equation}  
\label{macroVar}
V_c(t):=\frac{1}{N} \sum_{i=1}^N v_i(t).
\end{equation}
The microscopic variables are then taken as the fluctuations
around their mean values,
\begin{equation}  
\label{microVar}
w_i(t):=v_i(t) - V_c(t),
\qquad \mbox{for} \quad i=1,2,\dots,N.
\end{equation}
We denote $\w(t) = (w_1,w_2,\dots,w_N) \in \R^N$. Then we have 
\begin{equation}  
\w(t) = \vecv(t) - V_c(t) \, \e,
\qquad
\mbox{where} \quad \e := (1,1,\dots,1)^T\in\R^N.
\label{edefin}
\end{equation}
Since $\e$ is the eigenvector of $\A$ corresponding to the zero eigenvalue, 
we have $\A \w = \A \vecv$. Then (\ref{CS3}) can be rewritten as follows
\begin{equation}  
\d w_i 
= 
-\frac{\lambda}{N} (\widetilde\A\widetilde\w)_i\d t 
+ \frac{\sigma_i}{N} (\widetilde\A\widetilde\w)_i\d B_i^t 
- \d V_c.   
\label{microaux}
\end{equation}
The macroscopic variable $V_c$ satisfies the following lemma.

\medskip

\begin{lemma}\label{lem:macro}
Let $(\vecv(t))_{t\geq 0}$ be a solution of~$(\ref{CS0})$
subject to the deterministic constant initial 
datum~$(\ref{IC})$. Then $\,\E[V_c(t)] \equiv V_c(0)$ for 
$\,t\geq 0$ and $\,\E\left[ \int_{-\tau}^T |V_c(t)|^2 \d t \right] <\infty$
for all  $T<\infty$.
\end{lemma}

\medskip

\begproof
The boundedness of $\E\left[ \int_{-\tau}^T |V_c(t)|^2 \d t \right]$
follows directly from the definition (\ref{macroVar}) and the
martingale property of $\vecv(t)$ provided by Theorem \ref{thm:ex}.
Using (\ref{Psi}), we have
$$
\sum_{i=1}^N 
(\widetilde\A\widetilde\vecv)_i
= \e^T \widetilde\A\widetilde\vecv
= 0.
$$
Summing equations (\ref{CS3}), $i=1,2,\dots,N$,
using (\ref{macroVar}) and $\A \w = \A \vecv$, we obtain that 
the macroscopic dynamics is governed by the system
\begin{equation} 
\label{Vcdyn} 
\d V_c = 
\frac{1}{N^2} 
\sum_{i=1}^N \sigma_i (\widetilde\A\widetilde\w)_i \d B_i^t.
\end{equation}
After integration in time this implies
\[
\E[V_c(t)] = 
V_c(0) + \frac{1}{N^2} \sum_{i=1}^N \sigma_i \, 
\E\!\left[ \int_0^t (\A(s-\tau)\w(s-\tau))_i \d B_i^s \right]\qquad
   \mbox{for } t\geq 0.
\]
Since $f(s):=(\A(s-\tau)\w(s-\tau))_i$ is a martingale, we 
have $\E\left[ \int_0^t f(s) \d B_i^s \right] = 0$
(see Theorem 5.8 on p. 22 of \cite{Mao}). Thus we 
obtain $\E[V_c(t)] \equiv V_c(0)$.
\endproof
\medskip

\begin{remark}
Note that \eqref{microaux} and \eqref{Vcdyn} are expressed
in terms of the $\w$-variables only and so they form a
closed system, which is equivalent to \eqref{CS3}.
\end{remark}
\medskip

Clearly, due to~(\ref{microVar}), we have 
$\w \cdot \e = \sum_{i=1}^N w_i \equiv 0$.
Consequently, it is natural to introduce the decomposition
$\R^N = \langle\e\rangle \oplus \langle\e\rangle^\perp$, where $\e$
is given by~(\ref{edefin}).
We then have $\w(t)\in \langle\e\rangle^\perp$ for all $t\geq 0$.

\medskip

\begin{lemma}\label{lem:maxEig}
Let $\A\in\R^{N\times N}$, $N\geq 2$, be the matrix defined in $(\ref{Psi})$
and assume that $(\ref{PsiBoundSymm})$ and $(\ref{fiedlervect})$ hold. Then:

\smallskip

\noindent
{\rm (a)}
The maximal eigenvalue of $\A$ is bounded by $2\,(N-1)$.

\smallskip

\noindent
{\rm (b)}
We have $|\A\u|^2 \leq 2\,(N-1)\,\u^T\A\u$ for any vector $\u\in \R^N.$

\smallskip

\noindent
{\rm (c)}
We have
$\ell \,|\w|^2 \leq \w^T\A\w \leq  2\,(N-1)\,|\w|^2$
for any vector $\w\in \langle\e\rangle^\perp$.

\smallskip

\noindent
{\rm (d)}
For any vectors $\u$, $\w\in \langle\e\rangle^\perp$ and $\delta>0$ we have
\[
\u^T\A\w \leq  \frac{1}{2\delta}\,\u^T\A\u + \frac{\delta}{2}\,\w^T\A\w.
\]
\end{lemma}

\medskip

\noindent
\begproof
\begin{enumerate}[(a)]
\item The claim follows from the Gershgorin circle theorem.
Indeed, since $0 < \psi_{ij} \leq 1$, the diagonal entries 
satisfy $0 \leq \A_{ii} \leq N-1$,
and $\sum_{j\neq i} |\A_{ij}| = \A_{ii}$ for all $i=1,2,\dots,N$.
\item The smallest eigenvalue of $\A$ is zero with the corresponding 
eigenvector $\e$. The second smallest eigenvalue $\mu_2$ (Fiedler 
number) is assumed to be positive by (\ref{fiedlervect}).
Thus, $\A$ is a symmetric, positive operator on the space 
$\langle\e\rangle^\perp$ and there exists an orthonormal basis 
of $\langle\e\rangle^\perp$ composed of eigenvectors 
$\bxi_2,\bxi_3,\dots,\bxi_N$ of $\A$ corresponding to the 
positive eigenvalues $\mu_2,\mu_3,\dots,\mu_N$.
Then, every vector $\u \in \R^N$ can be decomposed as
\begin{equation}
\u 
= 
\frac{(\u \cdot \e) \, \e}{|\e|^2}
+
\sum_{i=2}^{N} (\u\cdot\bxi_i) \, \bxi_i.
\label{ubasis}
\end{equation}
Thus, due to the above bound on the eigenvalues 
$0 \leq \mu_i \leq 2(N-1)$, we have
$$
|\A\u|^2 
= 
\sum_{i=2}^N (\u\cdot\bxi_i)^2 \mu_i^2
\leq 2 \, (N-1) \sum_{i=2}^N (\u\cdot\bxi_i)^2 \mu_i
= 
2 \, (N-1) \, \u^T\A\u.
$$
\item
If $\w\in \langle\e\rangle^\perp$, then (\ref{ubasis}) implies
$$
\w 
= 
\sum_{i=2}^{N} (\w\cdot\bxi_i) \, \bxi_i,
\qquad \mbox{and} \qquad
\w^T\A\w 
= 
\sum_{i=2}^N (\w\cdot\bxi_i)^2\mu_i.
$$
Since nonzero eigenvalues are bounded from below by $\ell$ (using 
(\ref{fiedlervect})) and from above by part (a) of this lemma,
we obtain 
$$
\ell \, |\w|^2 = \ell \, \sum_{i=2}^N (\w\cdot\bxi_i)^2
\leq
\sum_{i=2}^N (\w\cdot\bxi_i)^2\mu_i 
\leq 2 \, (N-1) \, \sum_{i=2}^N (\w\cdot\bxi_i)^2 
= 2 \, (N-1) \, |\w|^2.
$$
\item
With the orthonormality of the basis $\bxi_2,\bxi_3,\dots,\bxi_N$ and 
the positivity of the eigenvalues $\mu_2,\mu_3,\dots,\mu_N$,
we have by the Cauchy-Schwartz inequality
\begin{eqnarray*}
\u^T\A\w 
&=& 
\left( \sum_{i=2}^{N} (\u\cdot\bxi_i) \, \bxi_i \right)^T \A 
\left(\sum_{i=2}^{N} (\w\cdot\bxi_i) \, \bxi_i \right)
= \sum_{i=2}^{N} (\u\cdot\bxi_i)(\w\cdot\bxi_i)\mu_i \\
&\leq& \left( \sum_{i=2}^{N} (\u\cdot\bxi_i)^2 
\mu_i \right)^\half \left( \sum_{i=2}^{N} 
(\w\cdot\bxi_i)^2 \mu_i \right)^\half
= (\u^T\A\u)^\half (\w^T\A\w)^\half,
\end{eqnarray*}
and with any $\delta>0$,
$$
(\u^T\A\u)^\half (\w^T\A\w)^\half 
\leq 
\frac{1}{2\delta} \, \u^T\A\u + \frac{\delta}{2} \, \w^T\A\w.
$$
\end{enumerate}
\vskip -5mm \rule{0pt}{0pt} \hfill \rule{0pt}{0pt}
\endproof

\subsection{Lyapunov functional}
\label{subsec:Lyapunov}
The proof of Theorem \ref{thm:flocking} relies on estimating the 
decay rate of the following Lyapunov functional for (\ref{microaux})--(\ref{Vcdyn}),
\begin{eqnarray}
\label{Lyapunov}
\qquad \!\mathscr{L}(t) := 
|\w(t)|^2 
&+& 
\frac{q}{N^2} \int_{t-\tau}^t \!\!\!|\A(s) \, \w(s)|^2 \d s 
\\
&+& 
\frac{p}{N^2} \int_{t-\tau}^t \int_\theta^t |\A(s-\tau) \, \w(s-\tau)|^2 \d s,  
\nonumber
\end{eqnarray}
where $p$, $q$ are positive constants depending on $\lambda$, 
$\tau$ and $\sigma_i$.

\medskip

\begin{lemma}\label{lem:Lyapunov}
Let the assumptions of Theorem $\ref{thm:flocking}$ be satisfied.
Then there exist positive constants $p$, $q$ and $\eps$ 
such that for every solution $(\w(t))_{t\geq 0}$ of
$(\ref{microaux})$--$(\ref{Vcdyn})$
the Lyapunov functional $(\ref{Lyapunov})$ satisfies
\begin{equation}   
\label{LyapDis}
\tot{}{t} \, \E[\mathscr{L}(t)] 
\, \leq \, -\frac{\eps}{N} \, \E \left[ \w^T\widetilde\A\w \right].
\end{equation}
\end{lemma}

\begproof
We apply the It\^{o} formula to calculate $\d w_i(t)^2$.
Note that the It\^{o} formula holds in its usual form 
also for systems of delayed SDE, see page 32 in \cite{Gillouzic}
and \cite{KM92, EN73, Mao2}.
Therefore, we obtain
\[
\d w_i(t)^2 
&=& -\frac{2\lambda}{N} w_i(\widetilde\A\widetilde\w)_i\d t 
+ \frac{2\sigma_i}{N} w_i(\widetilde\A\widetilde\w)_i\d B_i^t
- 2 w_i\d V_c \\
&& 
+ \frac{\sigma_i^2}{N^2} (\widetilde\A\widetilde w)_i^2 \d t
+ 
\frac{1}{N^4} \sum_{j=1}^N\sum_{k=1}^N 
\sigma_j \sigma_k (\widetilde\A\widetilde w)_j (\widetilde\A\widetilde w)_k 
\d B_j^t\d B_k^t \\
&& 
- 2\left( \frac{\sigma_i}{N} (\widetilde\A\widetilde w)_i \d B_i^t\right) 
\d V_c.
\]
With the identity $\d B^t_j\d B^t_k = \delta_{jk}\d t$ 
(formula (6.11) on p. 36 of \cite{Mao}), we have
\[
\frac{1}{N^4} \sum_{j=1}^N\sum_{k=1}^N \sigma_j \sigma_k 
(\widetilde\A\widetilde w)_j (\widetilde\A\widetilde w)_k \d B_j^t\d B_k^t
= 
\frac{1}{N^4} \sum_{j=1}^N \sigma^2_j (\widetilde\A\widetilde w)^2_j \d t.
\]
Consequently, summing over $i$, using $\w \cdot \e$ 
and the identity
\[
\sum_{i=1}^N \left( \frac{\sigma_i}{N} 
(\widetilde\A\widetilde w)_i \d B_i^t\right) \d V_c =
\frac{1}{N^3} \sum_{i=1}^N 
\sigma^2_i (\widetilde\A\widetilde w)^2_i \d t,
\]
we obtain
\[
\d |\w(t)|^2 &=& 
\left( -\frac{2\lambda}{N} \w^T \widetilde\A\widetilde \w 
+ \frac{N-1}{N^3}
\sum_{i=1}^N \sigma_i^2 (\widetilde\A\widetilde \w)_i^2 \right) \d t 
+ \frac{2}{N}\sum_{i=1}^N \sigma_i w_i (\widetilde\A\widetilde\w)_i 
\d B_i^t.
\]
Consequently, we have
\(
\d \mathscr{L}(t) &=&
\left( -\frac{2\lambda}{N} \w^T \widetilde\A\widetilde \w  + \frac{N-1}{N^3} 
\sum_{i=1}^N \sigma_i^2 (\widetilde\A\widetilde \w)_i^2 \right) \d t 
+ \frac{2}{N}\sum_{i=1}^N \sigma_i w_i (\widetilde\A\widetilde\w)_i 
\d B_i^t \nonumber \\
&& 
+ \frac{q}{N^2}\left( |\A\w|^2 - |\widetilde\A\widetilde\w|^2 \right)\d t
+ \frac{p}{N^2}\left(\tau |\widetilde\A\widetilde\w|^2 - \int_{t-\tau}^t 
|\A\w(s-\tau)|^2 \d s\right)\d t.
\label{d/dt EL}
\)
Our goal is to estimate $\tot{}{t}\E[\mathscr{L}(t)]$ from above.
First of all, we note that by the elementary property of the 
It\^{o} integral (Theorem 5.8 on p. 22 of \cite{Mao}),
\[
\E \left[ \frac{1}{N}\sum_{i=1}^N \sigma_i w_i (\widetilde\A\widetilde\w)_i \d B_i^t
\right] = 0.
\]
For the first term of the right-hand side in \eqref{d/dt EL}, we write
\[
-\frac{2 \, \lambda}{N} \w^T \widetilde\A\widetilde \w 
= -\frac{2\, \lambda}{N} \w^T \widetilde\A \w 
+ \frac{2\, \lambda}{N} \w^T \widetilde\A(\w-\widetilde \w)
\]
and apply Lemma \ref{lem:maxEig}(d) with $\delta>0$,
\[
\frac{2\, \lambda}{N} \w^T \widetilde\A(\w-\widetilde \w) 
\leq 
\frac{\lambda\, \delta^{-1}}{N} \w^T \widetilde\A \w 
+ \frac{\lambda \, \delta}{N} 
(\w-\widetilde \w)^T \widetilde\A (\w-\widetilde \w).
\]
Using Lemma \ref{lem:maxEig}(c), we have
\(   
\label{midstep}
\frac{\lambda \, \delta}{N} (\w-\widetilde \w)^T \widetilde\A (\w-\widetilde \w) 
\leq 2 \, \lambda \, \delta \, \frac{N-1}{N} \, |\w-\widetilde \w|^2.
\)
Now we write for $\w-\widetilde \w$, componentwise,
using \eqref{microaux},
\[
w_i-\widetilde w_i &=& \int_{t-\tau}^t \d w_i(s) \\
&=& - \frac{\lambda}{N}\int_{t-\tau}^t (\A(s-\tau) w(s-\tau))_i \d s
+ \frac{\sigma_i}{N} \int_{t-\tau}^t (\A(s-\tau) w(s-\tau))_i \d B_i^s\\
&&- \int_{t-\tau}^t \d V_c(s).
\]
Thus, we have for the expectation of the square
\[
\E |w_i-\widetilde w_i|^2 
&\leq& 3 \,
\E \!\left[
\frac{\lambda}{N}\int_{t-\tau}^t (\A(s-\tau) \w(s-\tau))_i \d s 
\right]^2
+ 3 \, 
\E \!\left[ 
\frac{\sigma_i}{N} \int_{t-\tau}^t (\A(s-\tau) \w(s-\tau))_i \d B_i^s 
\right]^2 \\
&& + 3 \, \E \!\left[ \int_{t-\tau}^t \d V_c(s) \right]^2.
\]
An application of the Cauchy-Schwartz inequality and Fubini's theorem for 
the first term of the right-hand side yields
\[
3 \, 
\E \!\left[
\frac{\lambda}{N}\int_{t-\tau}^t (\A(s-\tau) \w(s-\tau))_i \d s 
\right]^2  
\leq
\frac{3 \, \lambda^2}{N^2} \, \tau 
\int_{t-\tau}^t \E \!\left[ |(\A(s-\tau) \w(s-\tau))_i|^2 \right] \d s.
\]
For the second term we use the fundamental property of the It\^{o} 
integral (Theorem 5.8 on p. 22 of \cite{Mao}),
\[
3 \, \E 
\!\left[ 
\frac{\sigma_i}{N} \int_{t-\tau}^t (\A(s-\tau) \w(s-\tau))_i \d B_i^s 
\right]^2
= 
\frac{3 \, \sigma_i^2}{N^2} \int_{t-\tau}^t  
\E \!\left[ |(\A(s-\tau) \w(s-\tau))_i|^2 \right] \d s.
\]
Similarly, the third term is estimated as
\[
  3 \, \E \!\left[ \int_{t-\tau}^t \d V_c(s) \right]^2 &=&
   \frac{3}{N^4} \, \E \!\left[ \sum_{j=1}^N \int_{t-\tau}^t \sigma_j
     (\A(s-\tau) \w(s-\tau))_j \d B_j^s \right]^2 \\
   &\leq& \frac{3}{N^3} \sum_{j=1}^N \E \!\left[ \int_{t-\tau}^t \sigma_j
     (\A(s-\tau) \w(s-\tau))_j \d B_j^s \right]^2 \\
   &\leq& \frac{3\sigma_\mathrm{max}^2}{N^3} \sum_{j=1}^N \int_{t-\tau}^t
     \E \!\left[ |(\A(s-\tau) \w(s-\tau))_j|^2 \right] \d s.
\]
Thus, we get from \eqref{midstep}, estimating $\frac{N-1}{N}\leq 1$,
\[
\frac{\lambda \, \delta}{N} \, 
\E \!\left[ 
(\w-\widetilde \w)^T \widetilde\A (\w-\widetilde \w) 
\right]
\, \leq \, \frac{6 \, \lambda \, \delta}{N^2} 
 \left( \lambda^2 \tau + 2\sigma_\mathrm{max}^2 \right)
 \sum_{i=1}^N \int_{t-\tau}^t  
\E \!\left[ |(\A(s-\tau)\w(s-\tau))_i|^2 \right] \d s.
\]
An application of Lemma \ref{lem:maxEig}(b) gives
\[
\frac{q}{N^2} \, \E[ |\A\w|^2 ] 
\leq 
\frac{2 \, q \,(N-1)}{N^2} \, \E[\w^T\A\w] 
\leq 
\frac{2 \, q}{N} \, \E[\w^T\A\w].
\]
To balance this term with $-\frac{2\lambda}{N}\w^T\widetilde\A\w$, we use assumption $(\ref{Lipschitz})$
and Lemma \ref{lem:maxEig}(c) in
\[
   \w^T\A\w &=& \w^T(\A-\widetilde\A)\w + \w^T\widetilde\A\w \\
     &\leq& L\tau\w^T\w + \w^T\widetilde\A\w
     \leq \left(L\ell^{-1}\tau + 1\right) \w^T\widetilde\A\w.
\]

Collecting all the terms in \eqref{d/dt EL} finally 
leads to
\[
\tot{}{t} \, \E[\mathscr{L}(t)] &\leq& \frac{1}{N} 
\left[ -2 \, \lambda + \lambda \, \delta^{-1} + 2 \, q\left(L\ell^{-1}\tau + 1\right) \right] 
\E[\w^T \widetilde\A \w] 
+ \frac{1}{N^2} 
\left( \sigma_\mathrm{max}^2 + p \, \tau - q \right) 
\E |\widetilde\A\widetilde \w|^2 
\nonumber
\\
&+& \frac{1}{N^2} 
\left( 
6 \, \lambda \, \delta \,
(\lambda^2\tau + 2\sigma_\mathrm{max}^2) - p
\right)
\int_{t-\tau}^t  \E |\A(s-\tau)\w(s-\tau)|^2 \d s.
\]
We set
\begin{equation}
p = 6 \, \lambda \, \delta \, (\lambda^2\tau + 2\sigma_\mathrm{max}^2),
\qquad
q = \sigma_\mathrm{max}^2 + p \, \tau,
\label{definitionpq}
\end{equation}
then the above expression simplifies to
\begin{equation}
\tot{}{t} \, \E[\mathscr{L}(t)] \, \leq \, \frac{1}{N} 
\left[ -2 \, \lambda + \lambda \, \delta^{-1} + 2 \, q\left(L\ell^{-1}\tau + 1\right) \right]
\E[\w^T \widetilde\A \w].
\label{auxdleq2}
\end{equation}
We want $-2 \, \lambda + \lambda \, \delta^{-1} + 2 \, q\left(L\ell^{-1}\tau + 1\right) < 0$.
Substituting~(\ref{definitionpq}) into this inequality leads
to a third order polynomial inequality in $\tau$.
This polynomial has all positive coefficients but the 
zero order one, which is $c_0 := 2\sigma_\mathrm{max}^2 + \delta^{-1}\lambda - 2\lambda$.
If $(\ref{flockingCond_0})$ is satisfied, then choosing $\delta>0$ such that
\[
   \delta^{-1} < 2\lambda^{-1}(\lambda-\sigma_\mathrm{max}^2)
\]
makes $c_0$ negative. Consequently, there exists a $\tau_c >0$
such that for any $0 \leq \tau < \tau_c$,
\[
    -2 \, \lambda + \lambda \, \delta^{-1} + 2 \, q\left(L\ell^{-1}\tau + 1\right) = -\eps < 0.
\]
This completes the proof of $(\ref{LyapDis})$.

It remains to study the case when $\A$ is a constant matrix, i.e.
$L=0$ in $(\ref{Lipschitz})$.
Then $(\ref{auxdleq2})$ simplifies to
\[
\tot{}{t} \, \E[\mathscr{L}(t)] \, \leq \, \frac{1}{N} 
\left( -2 \, \lambda + \lambda \, \delta^{-1} + 2 \, q \right)
\E[\w^T \widetilde\A \w]
\]
and we have to find $\tau$ such that $-2 \, \lambda + \lambda \, \delta^{-1} + 2 \, q < 0$.
Again, substituting~(\ref{definitionpq}) for $p$ and $q$ leads to
\begin{equation}
\tau 
< \frac{\lambda \, (2-\delta^{-1})
-
2\,\sigma_\mathrm{max}^2}{12\,\lambda\,\delta
\,(\lambda^2\tau + 2\sigma_\mathrm{max}^2)}.
\label{auxdcondtau}
\end{equation}
The maximum value of the right hand side is obtained for 
$\delta = \lambda (\lambda-\sigma_\mathrm{max}^2 )^{-1}$
which is positive because of the first inequality in~(\ref{flockingCond}).
Substituting $\delta = \lambda (\lambda-\sigma_\mathrm{max}^2 )^{-1}$ 
into~(\ref{auxdcondtau}), we obtain
\[
\tau 
< 
\frac{(\lambda-\sigma_\mathrm{max}^2)^2}{12 \, \lambda^2 \,
(\lambda^2\tau + \sigma_\mathrm{max}^2)}.
\]
Finally, resolving in $\tau$ leads to
\[
\tau < \frac{1}{\lambda^2}
\left( - \, \sigma_\mathrm{max}^2 
+ 
\sqrt{\, \sigma_\mathrm{max}^4 
+ \frac{1}{12} \, (\lambda-\sigma_\mathrm{max}^2)^2} \right).
\]
If the above sharp inequality is satisfied, there exists an $\eps>0$ such that
$-\eps = -2 \, \lambda + \lambda \, \delta^{-1} + 2 \, q$ and, consequently,
(\ref{LyapDis}) holds.
\endproof

\subsection{Proof of Theorem \ref{thm:flocking}}
An integration of \eqref{LyapDis} in time gives
\[
\E[|\w|^2(t)] 
&\leq& 
\E[\mathscr{L}(t)] 
= \E[\mathscr{L}(0)] + \int_0^t \tot{}{t} \, \E[\mathscr{L}(s)] \d s \\
&\leq& \E[\mathscr{L}(0)] -\frac{\eps}{N} \int_0^t \E[\w^T(s)\A(s)\w(s)] \d s.
\]
An application of Lemma \ref{lem:maxEig}(c) gives then
\[
\E[|\w|^2(t)] \leq
 \E[\mathscr{L}(0)] 
-
\frac{\eps \, \ell}{N} \int_0^t \E[|\w(s)|^2] \d s,
\]
so that the last integral is convergent as $t\to\infty$ and, 
consequently, $\lim_{t\to\infty} \E[w_i(t)] = 0$ for 
all $i=1,2,\dots,N$. Using \eqref{microVar}, we obtain
\[
\lim_{t\to\infty} \left| \E[v_i(t)]- \E[v_j(t)] \right| 
= 
\lim_{t\to\infty} \left| \E[w_i(t)]- \E[w_j(t)] \right| 
= 
0,
\]
and we conclude that asymptotic flocking in the sense of 
Definition \ref{def:flocking} takes place.
\vbox{\hrule height0.6pt\hbox{%
   \vrule height1.3ex width0.6pt\hskip0.8ex
   \vrule width0.6pt}\hrule height0.6pt
}

\medskip

\begin{remark}
The original definition of asymptotic flocking~{\rm \cite{CS1,CS2}} involves,
in addition to our Definition $\ref{def:flocking}$, also the 
group formation property for~$(\ref{CS01})$--$(\ref{CS02})$ given by
\[
\sup_{t\geq 0} \left| \E[x_i(t)]- \E[X_c(t)] \right| 
< 
\infty,
\qquad\qquad
\mbox{for all }
\quad 
i=1,2,\dots,N,
\]
where $X_c(t):=\frac{1}{N}\sum_{i=1}^N x_i(t)$.
The standard way~{\rm \cite{CS1, CS2, Tadmor-Ha, Ha-Liu}} of proving this 
result would be to estimate
\[
\left| \E[x_i(t)]- \E[X_c(t)] \right| 
\leq \left|x_i^0 - X_c(0)\right| + \int_0^t \E[|w_i(s)|] \d s
\]
and employ a bootstrapping argument to show that 
$\int_0^\infty \E[|w_i(s)|] \d s < \infty$.
However, as noted above, it is not clear how to apply the 
bootstrapping argument in our setting. Note that we have 
$\int_0^\infty \E[|w_i(s)|^2] \d s < \infty$,
but it does not imply $\int_0^\infty \E[|w_i(s)|] \d s < \infty$.
\end{remark}
\medskip

\begin{remark}
Let us note that the flocking conditions $(\ref{flockingCond_0})$ and $(\ref{flockingCond})$, 
expressed in terms of $\lambda$ and $\sigma_\mathrm{max}$,
translate into
\[
\lambda \, \kappa_\mathrm{max}^2 <1,\qquad
\tau < -\kappa_\mathrm{max}^2 
+ \sqrt{\kappa_\mathrm{max}^4 
+ \frac{1}{12\lambda^2} \, (1-\lambda \, \kappa_\mathrm{max}^2)^2}
\]
in terms of $\lambda$ and 
$\kappa_\mathrm{max} := \sigma_\mathrm{max}/\lambda$,
where $\kappa_\mathrm{max}$ can also be written as
$
\kappa_\mathrm{max}
=
\max \{\kappa_1,\kappa_2,\dots,\kappa_N\}
$ 
for $\kappa_i$, $i=1,2,\dots,N$, defined by $(\ref{omega_i})$.
\end{remark}

\medskip

\subsection{Application to asymptotic behaviour of delayed 
geometric Brownian motion}\label{subsec:delayedGBM}
Our analysis provides information about the asymptotic behaviour of
the delayed geometric Brownian motion \eqref{simplified}
which is, to our best knowledge, new.
We just modify the proof of Lemma \ref{lem:Lyapunov}
with the obvious simplifications due to the fact that $\A(t)\equiv 1$.
This leads to a slight improvement in the flocking condition.

\medskip

\begin{lemma}\label{lem:delayedGBM}
Let the parameters 
$\lambda>0$, $\sigma \in \R$ and $\tau\geq 0$ satisfy
\(   
\label{condDGBM}
\sigma^2 < 2 \lambda,
\qquad 
\tau < \frac{1}{4 \, \lambda^2}
\left( -2 \, \sigma^2 + \sqrt{4 \, \sigma^4 + 2(2 \, \lambda-\sigma^2)^2} 
\right).
\)
Then the solutions of the delayed geometric Brownian motion equation
$(\ref{simplified})$ satisfy
\[
    \lim_{t\to\infty} \E[|w|^2] = 0.
\]
\end{lemma}

Let us note that the above result is suboptimal for the 
deterministic case.
Indeed, setting $\sigma:=0$, \eqref{condDGBM} reduces 
to $\lambda\tau < \sqrt{2}/2$.
However, it is known \cite{Smith} that solutions of the 
delayed ODE $\dot w = -\lambda\widetilde w$
asymptotically converge to zero if $\lambda\tau < \pi/2$.
On the other hand, if there is no delay, i.e., $\tau = 0$,
the condition \eqref{condDGBM} reduces to $\sigma^2 < 2\lambda$,
which is the sharp condition for asymptotic vanishing of 
mean square fluctuations of geometric Brownian motion.

\section{Numerical experiments}\label{sec:Numerics}
We provide results of numerical experiments for the models
considered in this paper with focus on their asymptotic 
behaviour. First, we illustrate that one has to be cautious
when interpreting the numerical results as indications about 
the ``true'' asymptotic behaviour of the solution, because
implementations of Monte-Carlo algorithms for geometric 
Brownian motion lead to systematic underestimation
of the moments of the true solution, see Section \ref{subsec:MCdefect}.
Keeping this systematic defect in mind, we will resort
to weak methods for simulation of our SDEs
and study their \emph{numerical} asymptotic behaviour.
In Section \ref{subsec:numDGBM}, we resort to the 
delayed geometric Brownian motion \eqref{simplified},
which can be seen as a toy model of \eqref{CS01}--\eqref{CS02},
and find combinations of parameter values that guarantee
\emph{numerical} asymptotic decay of the solution.
In Section \ref{subsec:numCS0}, we then perform
numerical simulations of the velocity alignment system \eqref{CS0}
with fixed communication rates, and, finally,
in Section \ref{subsec:numCS} we focus
on the full system \eqref{CS01}--\eqref{CS02}.

\subsection{Analysis of the Monte-Carlo method for geometric 
Brownian motion}\label{subsec:MCdefect}
In this section we estimate the systematic error produced
by numerical implementations of the Monte-Carlo algorithm
for geometric Brownian motion without delay. We show that
computer simulations underestimate the mean square fluctuations 
of the process due to the fact that the numerical implementation 
does not capture large deviations (extreme outliers), and the error
grows exponentially in time.

\begin{figure}
\vskip 3mm
\centerline{
\hskip 3mm
\includegraphics[width=0.48\columnwidth]{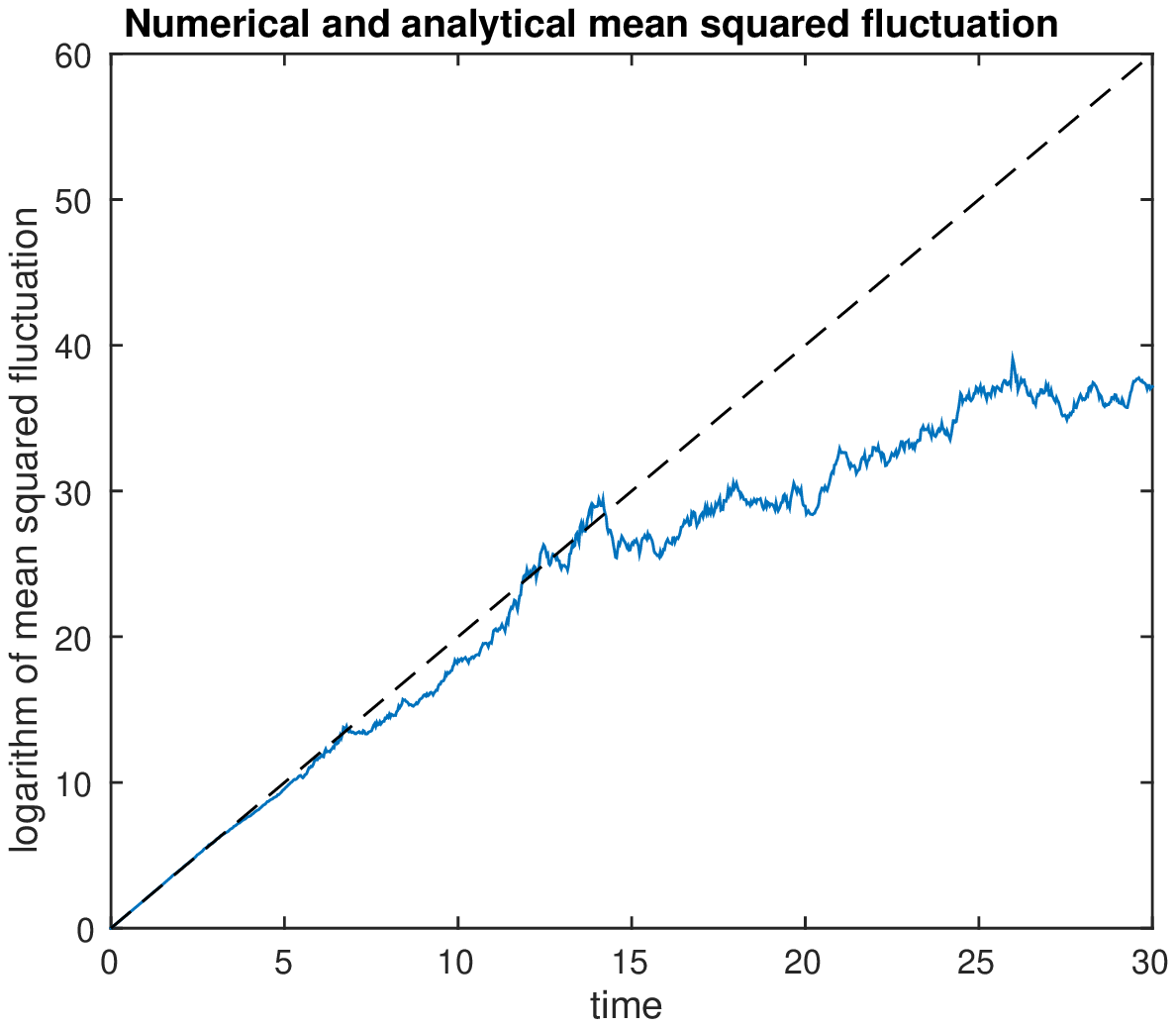}
\hskip 3mm
\includegraphics[width=0.48\columnwidth]{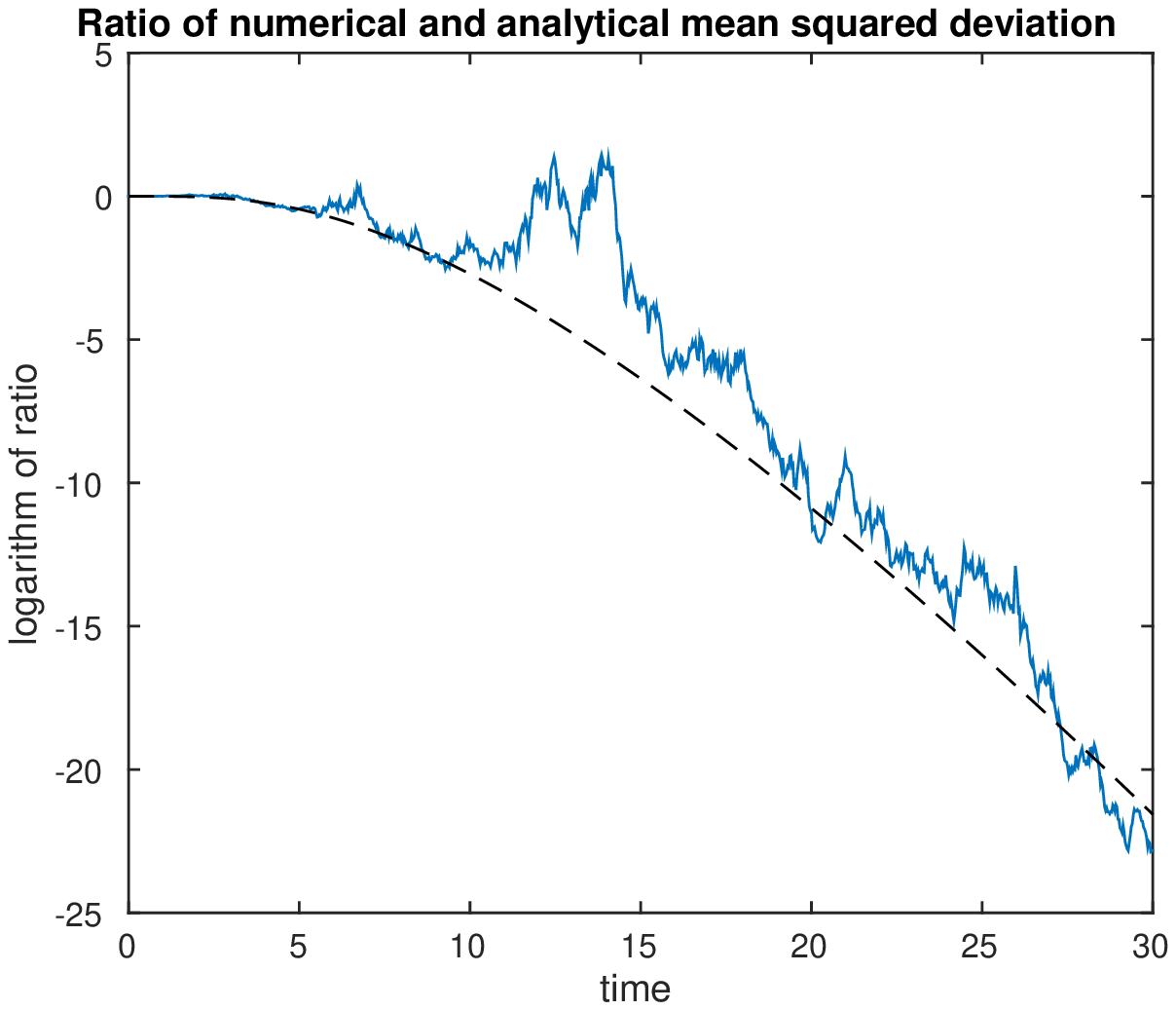}
}
\vskip -5.3cm
\leftline{\hskip 2mm (a) \hskip 6.15cm (b)}
\vskip 4.7cm
\caption{{\rm (a)} {\it Logarithm of the simulated mean squared 
fluctuations $\log(\E[\bar v^2(t)])$ (solid line),
and the analytical result $\log(\E[v^2(t)])=4 \, \lambda \, t$ 
(dashed line).}
{\rm (b)} {\it Logarithm of the ratio of simulated and analytically 
calculated mean squared fluctuations
$\E[\bar v^2(t)]/\E[v^2(t)]$ (solid line), and the 
theoretically calculated curve (dashed line) given by
the right-hand side of $(\ref{ratio})$.
The Monte-Carlo simulation for $z(t)$ was performed with $10^6$ paths
of the process $(\ref{xi})$ with $z(0)=0$, $\lambda=0.5$ and $\sigma=1$
on the time interval $[0,30]$ divided into $10^3$ equidistant 
sampling points.}
\label{fig:MSfluct}}
\end{figure}

Let us consider the one-dimensional Brownian motion with drift,
\(   
\label{xi}
\d z  = (\lambda-\sigma^2/2)\d t + \sigma \d B^t,
\qquad\qquad 
z(0) = 0.
\)
Then, defining $v(t):=\exp[z(t)]$, we have by the It\^{o} formula
\(   
\label{v_t_SDE}
\d v = \lambda \, v \d t + \sigma \, v \d B^t,
\qquad \qquad
v(0) = 1.
\)
For simplicity, we perform a model calculation with 
$2 \, \lambda=\sigma^2 > 0$,
so that
\( \label{xi_and_v}
\d z
= 
\sqrt{2 \, \lambda} \d B^t, \qquad \qquad  
\d v = \lambda \, v \d t +\sqrt{2 \, \lambda} \, v \d B^t.
\)
Then, the density $u(t,x)$ of the process $z$ is given by
\[
u(t,x) 
= 
\frac{1}{\sqrt{4  \pi  \lambda \, t}} 
\, 
\exp \!\left( -\frac{x^2}{4 \lambda  t} \right).
\]
Moreover, we have
\[
\E[z(t)] = 0, \qquad 
\E[z^2(t)] = 2 \, \lambda \, t, \qquad
\E[v(t)] = \exp(\lambda \, t), \qquad 
\E[v^2(t)] = \exp( 4 \, \lambda \, t).
\]
We assume that our numerical scheme produces 
approximations $\bar{z}$ of the process $z$
that excludes the extreme outliers, i.e., 
$\mbox{Prob}(|\bar{z}(t)|>\alpha(t)) = 0$ for some
$\alpha=\alpha(t)$. In particular, we consider a properly 
scaled cut-off of the density $u(t,x)$ such that the 
probability of the extreme outliers $\mbox{Prob}(|z(t)|>\alpha(t))$ 
remains constant in time. This leads to 
$\alpha(t) = \eta \, \sqrt{4 \, \lambda \, t}$ for some $\eta >0$, 
since
\[
\mbox{Prob}(|z(t)|>\eta \, \sqrt{4 \, \lambda \, t}) 
= 
2 
\int_{\eta \, 
\sqrt{4 \, \lambda \, t}}^{\infty} 
u(t,x) \d x = \mbox{erfc}(\eta),
\]
where $\mbox{erfc}(\eta)= \frac{2}{\sqrt{\pi}} 
\int_\eta^\infty \exp \!\left( - x^2 \right) \d x$
is the complementary error function.
Consequently, we turn $u(t,x)$ into the truncated 
probability density
\[
\bar u_\eta(t,x) 
:= 
\frac{\chi_{[-\alpha(t), \alpha(t)]}(x)}{\mbox{erf}(\eta)} \, u(t,x),
\]
where $\chi_{[-\alpha(t), \alpha(t)]}$ 
is the characteristic function of interval
$[-\alpha(t), \alpha(t)]$ and
$\mbox{erf}(\eta) = 1 - \mbox{erfc}(\eta)$ is the error function.
Let us denote by $\bar{z}(t)$ the process with the density 
$\bar u_\eta(t,x)$ for a fixed $\eta>0$,
and $\bar v(t) := \exp(\bar{z}(t))$.
A simple calculation then reveals that
\[
\E[\bar v^2(t)] 
= \int_{-\infty}^{\infty} \exp(2x) \d\bar 
u_{\eta}(x)
= 
\frac{\exp(4 \, \lambda \,t)}{2 \, \mbox{erf}(\eta)} 
\left( \mbox{erfc}\!\left(-\eta-\sqrt{4\lambda t}\right) 
- \mbox{erfc}\!\left(\eta-\sqrt{4\lambda t}\right) \right).
\]
Consequently, since $\E[v^2(t)] = \exp(4 \, \lambda \, t)$, 
the numerical method produces the relative error
\(   \label{ratio}
\frac{\E[\bar v^2(t)]}{\E[v^2(t)]}
= 
\frac{\mbox{erfc}(-\eta-\sqrt{4\lambda t}) 
- \mbox{erfc}(\eta-\sqrt{4\lambda t})}{2 \, \mbox{erf}(\eta)}.
\)
This ratio is equal to one for $t=0$. Using the mean value theorem,
we obtain the asymptotic behaviour of the ratio for large times,
\[
\frac{\E[\bar v^2(t)]}{\E[v^2(t)]}
\approx
\frac{2 \, \eta}{\sqrt{\pi} \, \mbox{erf}(\eta)} \exp(-4\lambda t),
\qquad \quad \mbox{as } \;t\to\infty.
\]
Consequently, any implementation of the Monte Carlo method
excluding large deviations will underestimate the true value 
of $\E[v^2(t)]$ by an exponentially growing factor in time.
Let us note that this is also true for any moment of $v$ 
and with general parameters $\lambda$ and $\sigma$.
We illustrate this fact using a numerical simulation.
We perform a Monte Carlo simulation in {\tt Matlab} with $10^6$ 
paths of the process \eqref{xi} on the time interval $[0,T]$
with $T=30$ and $10^3$ equidistant sampling points.
We impose the initial condition $z(0) = 0$ and the parameter 
values $\lambda=0.5$ and $\sigma=1$.
Consequently, $z(t)$ is the Wiener process $B^t$,
and for its numerical approximation $\bar{z}$ we use
the built-in {\tt Matlab} procedure {\tt normrnd}
that generates normally distributed random numbers.
We calculate $\bar v(t):=\exp(\bar{z}(t))$ and evaluate the
mean squared fluctuations $\E[\bar v^2(t)]$.
We plot its logarithm as the solid curve in 
Figure \ref{fig:MSfluct}(a), compared to the analytical
curve $\log(\E[v^2(t)])=4 \, \lambda \, t$ (dashed line). We observe 
the exponential in time divergence of the two curves.
This is well described by our formula (\ref{ratio}),
as illustrated in Figure \ref{fig:MSfluct}(b).
For the calculation of the cut-off parameter $\eta$
we use the maximal value attained by the actual numerical realization
of the stochastic process, i.e., we set
$\eta := \max_{t\in(0,T]} \frac{|\bar{z}(t)|}{\sqrt{4 \, \lambda \, t}}$.
We then plot the logarithm of the ratio $\E[\bar v^2(t)]/\E[v^2(t)]$
and the theoretically calculated curve given by
the right-hand side of (\ref{ratio}).
We observe a good match between the two curves.

This systematic discrepancy between the analytical formulas and 
Monte Carlo simulations originates in the heavy tailed distribution 
of the geometric Brownian motion, and is a well studied topic, see, 
e.g., the survey~\cite{Juneja}. Importance sampling and rare-event 
simulation techniques would be the methods of choice to overcome 
this problem, however, their implementation is beyond the scope of 
our paper. Instead, we argue that the motivation for studying the 
stochastic system~(\ref{CS01})--(\ref{CS02}) and 
its simplification~(\ref{simplified})
comes from the fact that they ought to represent models of some 
real (physical or biological) phenomena. In real life situations these
extreme events with exponentially low probabilities can be unphysical
and the presented Monte Carlo simulation can be considered a more 
appropriate description of reality than the SDEs. 
We adopt this point of view for the forthcoming numerical studies 
in Sections~\ref{subsec:numDGBM},~\ref{subsec:numCS0} and~\ref{subsec:numCS},
and accept the fact their results may not be directly related to 
statements of Theorem~\ref{thm:flocking} and Lemma~\ref{lem:delayedGBM}.

\subsection{Numerical study of delayed geometric Brownian motion}
\label{subsec:numDGBM}
\begin{figure}
\centerline{\resizebox*{0.76\linewidth}{!}%
{\includegraphics{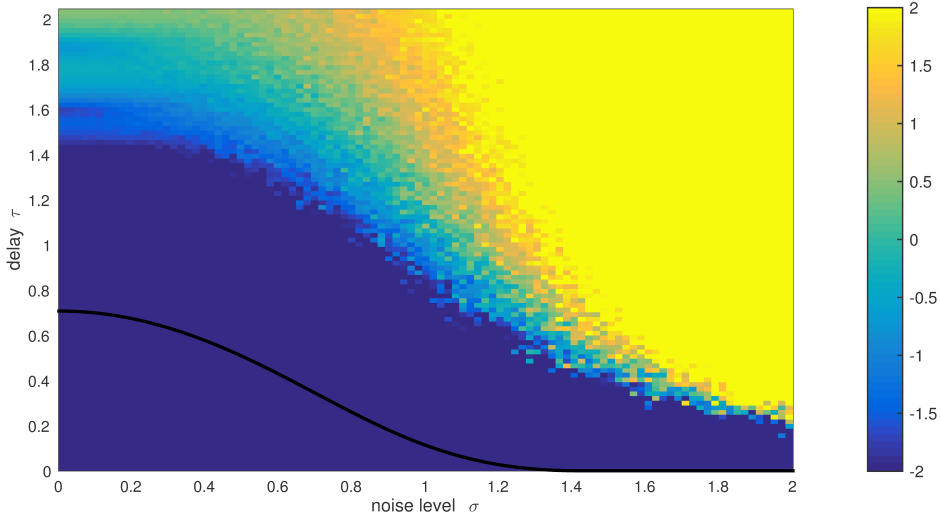}}}
\caption{Results of Monte Carlo simulations of the delayed 
SDE \eqref{simplified2}
with $Q=100$ paths for $(\sigma,\tau)\in[0,2]\times[0,2]$. 
The background colour encodes $\log(I_{\sigma,\tau})$. The 
region of ``numerical convergence'' is dark blue. The solid 
line indicates the critical values of $\tau$ given by 
formula \eqref{condDGBM} with $\lambda=1$ as a function of $\sigma$.
\label{fig:noise-delay}}
\end{figure}

Using $\lambda = 1$, the delayed SDE \eqref{simplified} can be 
equivalently written as
\begin{equation}
\label{simplified2}
\d w = - \widetilde w \, \d t + \sigma \, \widetilde w \, \d B^t,
\end{equation}
where $\sigma \ge 0$ and $\tau \ge 0$ are nonnegative parameters. 
We perform a systematic numerical study of the delayed SDE 
\eqref{simplified2} to characterize the asymptotic behaviour of 
its solutions in dependence on the values of the parameters 
$\sigma$ and $\tau$. In particular, we divide the domain 
$[0,2]\times[0,2]$ for $(\sigma,\tau)$ into $100\times 100$ 
equidistant $(\sigma,\tau)$-pairs. For each pair of the parameter 
values we perform a Monte Carlo simulation for \eqref{simplified2} 
with $Q=100$ paths over the time interval $[0,T]$ with $T=30$ 
and timestep $\Delta t = 10^{-3}$. We impose the constant 
deterministic initial condition $w(t) \equiv 1$ for
$t \in (-\tau,0]$. For discretization of \eqref{simplified2} 
we use the Euler-Maruyama method, i.e., the discrete scheme is
\(   \label{discrete}
w_{t_{k+1}} 
= 
w_{t_k} - \Delta t \, w_{t_k-\tau} + \sigma \, 
\sqrt{\Delta t} \, \mathcal{N}_{0,1},
\qquad k=1,2,\dots,K,
\)
subject to the initial condition $w_t \equiv 1$ for $t\leq 0$. 
Here $K=T/\Delta t$ denotes the total number of timesteps,
$t_k = k\Delta t$, and $\mathcal{N}_{0,1}$ a normally 
distributed random variable with zero mean and unit variance.
Note that the values of $\tau$ are chosen to be integer multiples 
of $\Delta t$, so that $t_k-\tau = t_l$ for some $l\in\mathbb{Z}$.
For each $(\sigma,\tau)$-pair and each path $q$ of the 
Monte Carlo simulation we calculate the ``indicator''
\[
I_{\sigma,\tau} := \frac{1}{Q}\sum_{q=1}^Q \left( 
\Delta t \sum_{k=(T-1)/\Delta t}^{T/\Delta t} |w^q_{t_k}|^2 \right)^{1/2},
\]
where $w^q_t$ is the $q$-th path in the Monte Carlo 
simulation of \eqref{discrete}.
The background colour in Figure~\ref{fig:noise-delay} encodes 
the logarithm of $I_{\sigma,\tau}$. To define a region 
of ``numerical convergence'', we choose a threshold $\Theta$
such that $I_{0,\tau_c} = \Theta$ for the delay $\tau_c = \frac{\pi}{2}$
that is critical for the problem without noise ($\sigma=0$).
In our case this led to $\Theta \simeq 10^{-2}$.
The region of ``numerical convergence'' is marked dark blue 
in Figure~\ref{fig:noise-delay}.
We observe the decrease of the critical value of the delay 
with increasing level of noise.
For comparison, the critical values of $\tau$ given by 
formula \eqref{condDGBM} as a function of $\sigma$
are indicated by the solid line.

\subsection{Numerical study of the system \eqref{CS0} with 
fixed communication matrix}
\label{subsec:numCS0}
We present results of numerical simulations of the system \eqref{CS0}
in the one-dimensional setting $d=1$, where we fix the communication rates
to $\psi_{ij}\equiv 1$ for all $i,j=1,2,\dots,N$,
i.e., every agent communicates with all others at the same rate.
Consequently, the communication matrix $\A$ has the off-diagonal
entries $\A_{ij}=-1$, $i\neq j$, and $\A_{ii}=N-1$.
It only has two eigenvalues, $0$ and $N$. Consequently,
its Fiedler number is $\mu_2=N$ and
we can choose $\ell:=N$ in $(\ref{fiedlervect})$.
In this setting, we can directly compare our analytical result,
Theorem \ref{thm:flocking}, with numerical simulations.

We will be considering even number of agents $N=2K$, 
in particular, $N\in\{2,20\}$,
and prescribe the initial datum
\(   \label{ICv}
v_i(t) \equiv 
\left\{
\begin{matrix}
-1 & \quad\mbox{for } i=1,2,\dots,K; \hfill \\
1 & \quad\mbox{for } i=K+1,K+2,\dots,N;
\end{matrix}
\right.
\)
for $t\in(-\tau,0]$. Although the asymptotic behaviour of the solutions 
in general depends on the particular choice of the initial datum,
a systematic study of this dependence is beyond the scope of this paper.
Therefore we only consider the ``generic'' choice of
initial conditions~(\ref{ICv}).

\begin{figure}
\vskip 3mm
\centerline{
\hskip 3mm
\includegraphics[width=0.48\columnwidth]{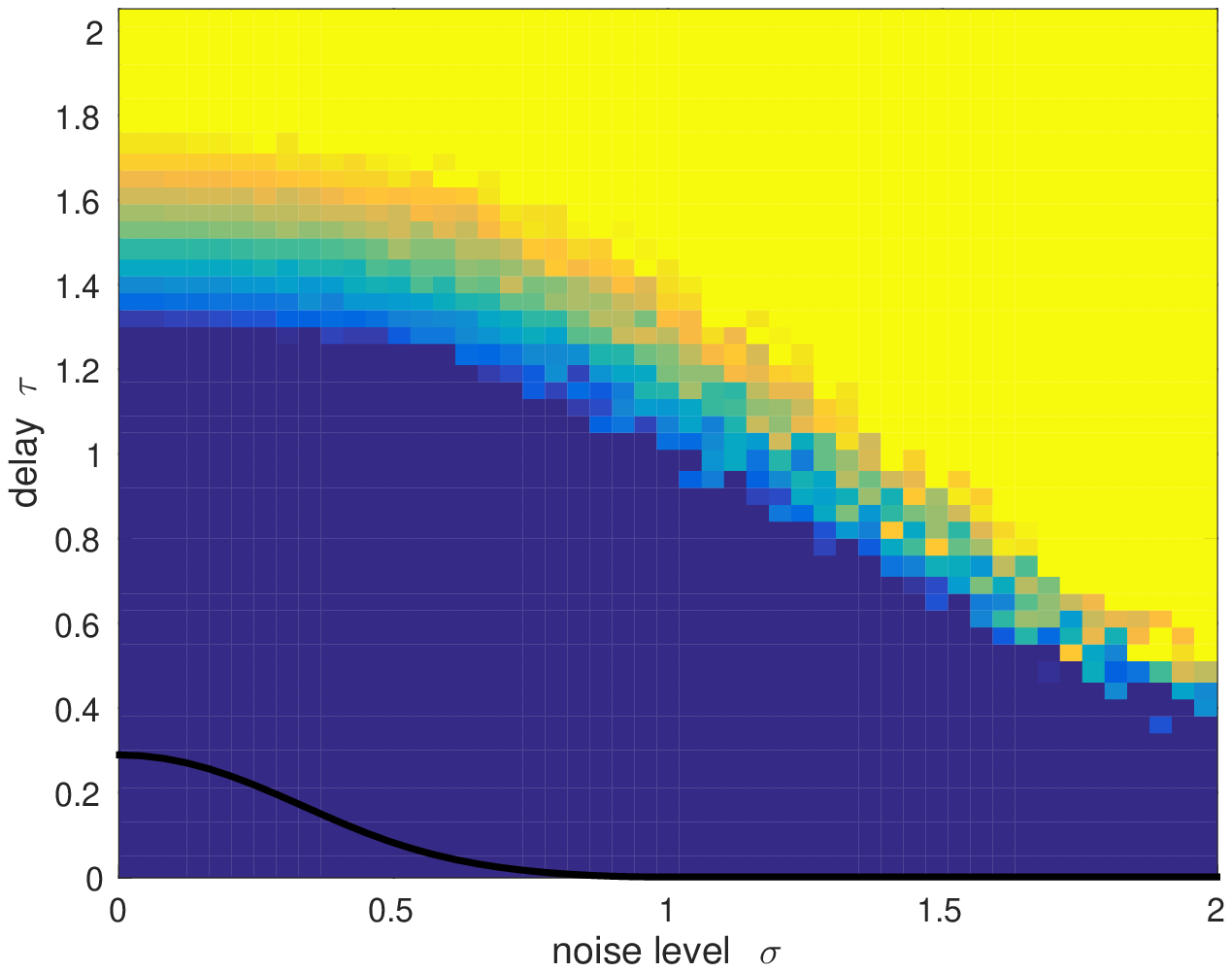}
\hskip 3mm
\includegraphics[width=0.48\columnwidth]{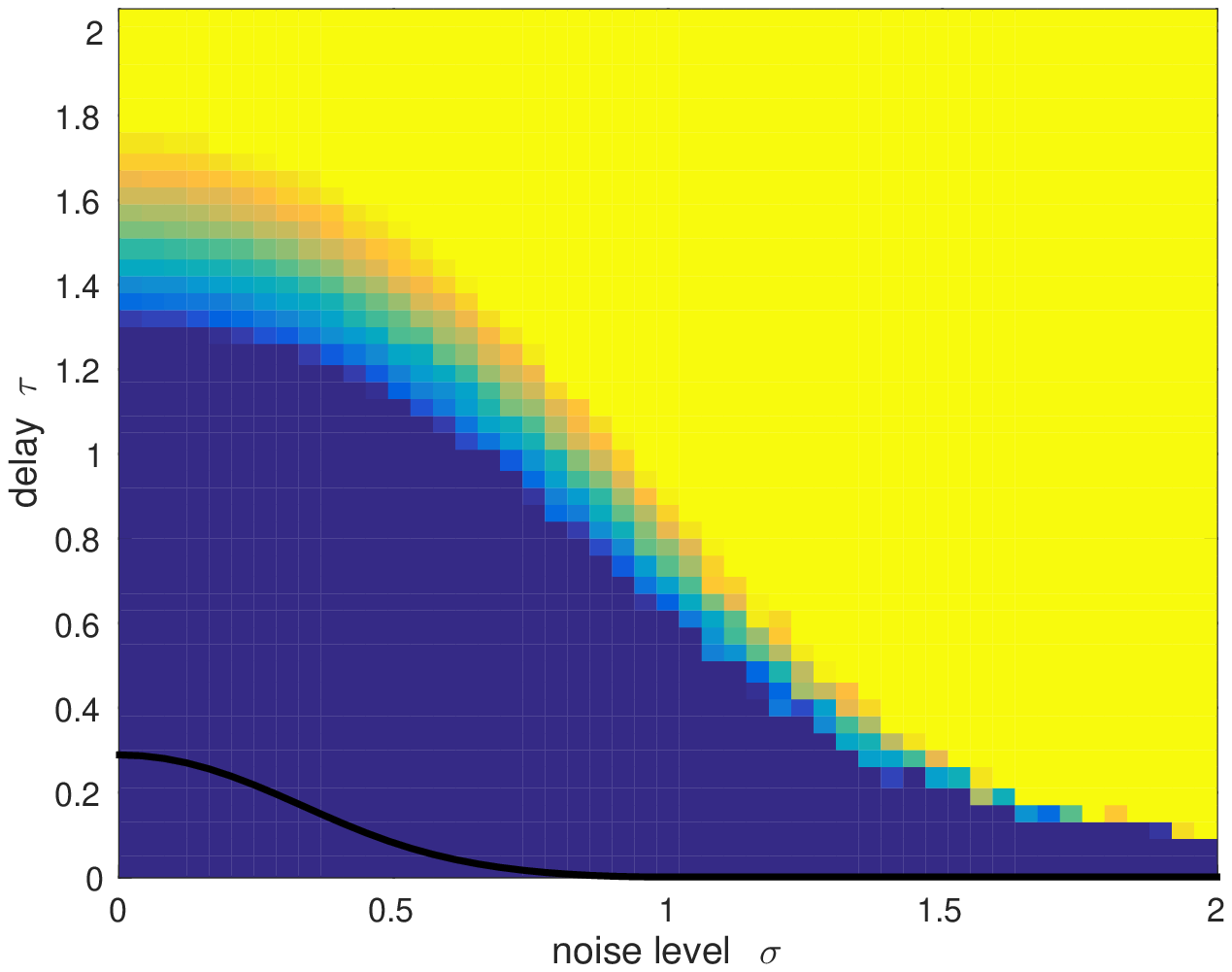}
}
\vskip -5cm
\leftline{\hskip 2mm (a) $N=2$ \hskip 5.15cm (b) $N=20$}
\vskip 4.5cm
\centerline{
\resizebox*{0.6\linewidth}{!}{\includegraphics{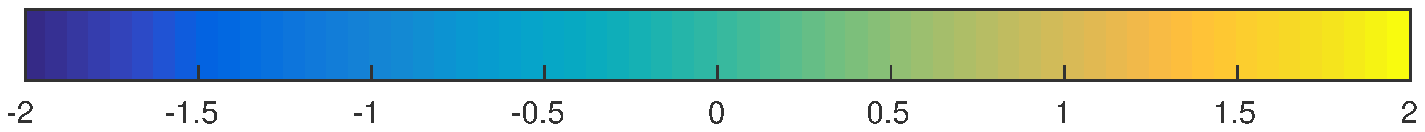}}
}
\caption{
{\it
Decadic logarithm of the indicator $I_{\sigma,\tau}$, given by
\eqref{indicator},
for Monte Carlo simulations of the system \eqref{CS0} with $Q=100$
paths on the time interval $[0,30]$ subject to the initial condition 
\eqref{ICv},
with $\lambda=1$, $(\sigma,\tau)\in[0,2]\times[0,2]$.
The dark blue regions (colour online) indicate ``numerical flocking''.
The solid line indicates the critical value $\tau_c$ given
by formula \eqref{flockingCond} for $\lambda=1$ as a function of $\sigma$.
The number of individuals is:} {\rm (a)} $N=2$;
{\it and} {\rm (b)} $N=20$.
\label{fig:CS0-syst}}
\end{figure}

We perform Monte Carlo simulations of the system \eqref{CS0}
with $N\in\{2,20\}$,  $\sigma_i=\sigma$ for all 
$i=1,2,\dots,N$ and $\lambda=1$
(other values of $\lambda$ can be achieved by rescaling of 
$\sigma$ and time).
We divide the domain $[0,2]\times[0,2]$ for $(\sigma,\tau)$
into $50\times 50$ equidistant $(\sigma,\tau)$-pairs. For each pair 
of the parameter values we perform a Monte Carlo simulation 
with $Q=100$ paths over the time interval
$[0,T]$ with $T=30$. We use the Euler-Maruyama method
for discretization of \eqref{CS0} with timestep $\Delta t = 10^{-3}$.
To classify the asymptotic behaviour of the solution, we again 
define the ``indicator''
\(   \label{indicator}
I_{\sigma,\tau} 
:= \frac{1}{Q} \sum_{q=1}^Q \left( \frac{\Delta t}{N} 
\sum_{k=(T-1)/\Delta t}^{T/\Delta t} |\vecv^q_{t_k}|^2 \right)^{1/2}
\)
where $\vecv^q_{t_k}$ is the $q$-th path in the Monte Carlo simulation 
of \eqref{CS01}--\eqref{CS02} at time $t_k=k\Delta t$.
We say that \emph{numerical flocking} takes place when 
$I_{\sigma,\tau} < 10^{-2}$.
The background colour in Figure~\ref{fig:CS0-syst} encodes 
the decadic logarithm of the indicator,
and the dark blue region indicates numerical flocking.
We observe that the region of numerical flocking is only 
weakly influenced by the number of agents $N$.
This is in agreement with the fact that the flocking 
condition \eqref{flockingCond} in Theorem \ref{thm:flocking}
does not depend on $N$.
The increased smoothness of the colour transition when $N=20$ 
is a consequence of the law of large numbers.
For comparison, the critical value $\tau_c$ given by 
formula \eqref{flockingCond} for $\lambda=1$ as a function of 
$\sigma$ is indicated by the solid line in both panels. 
The comparison with the numerical results suggests that
the condition \eqref{flockingCond} is far from optimal.

\subsection{Numerical study of the delayed Cucker-Smale system 
with multiplicative noise}

\label{subsec:numCS}

\begin{figure}
{\centering
\resizebox*{\linewidth}{!}{\includegraphics{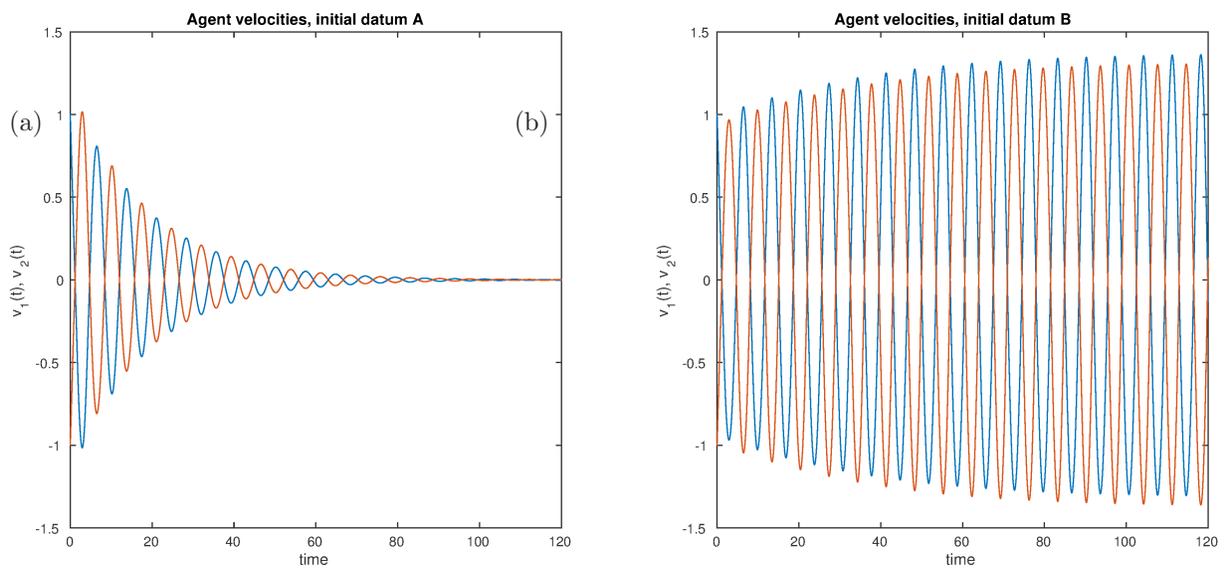}}
\par}
\vskip -6.6cm
\leftline{\hskip 2mm (a) \hskip 6.15cm (b)}
\vskip 6cm
\caption{{
\it Numerical simulations of the system 
\eqref{CS01}--\eqref{CS02}, \eqref{commRate}
with parameter values $N=2$, $\lambda=1$, $\beta=0.1$, 
$\sigma=0$ and $\tau = 1.75$.
Both simulations are performed on the time interval $[0,120]$ 
with discrete timestep $\Delta t = 10^{-3}$.
The initial condition for $v$ is $v_1(t)\equiv 1$, 
$v_2(t)\equiv -1$ for $t\in(-\tau,0]$
in both cases. The initial condition for $x$ is:}
{\rm (a)} $x_1(t)\equiv -1$, $x_2(t)\equiv 1$;
{\rm (b)} $x_1(t)\equiv 1$, $x_2(t)\equiv -1$.
The plots show the velocities of the two agents (red and  
blue, colour online) as functions of time.
\label{fig:IC}}
\end{figure}

Finally, we present results of numerical simulations of the 
system \eqref{CS01}--\eqref{CS02}
in the one-dimensional setting $d=1$ with the communication 
rates $\psi_{ij}=\psi(|x_i-x_j|)$
and $\psi$ given by \eqref{commRate}.
As in Section~\ref{subsec:numCS0}, our goal is to characterize 
the asymptotic behaviour of the solutions
in dependence on the parameter values, however, we are facing 
additional difficulties here.
In particular, the asymptotic behaviour of the solution may 
depend nontrivially on the initial condition,
as we show in Figure~\ref{fig:IC}. Since a systematic study 
taking this effect into account 
is beyond the scope of this paper, we will impose the same 
type of initial condition for all our simulations.
In particular, we prescribe constant zero value for 
the $\x$-variables,
\(   
\label{ICx}
x_i(t) \equiv 0 \quad\mbox{for } t\in(-\tau,0],\; i=1,2,\dots,N.
\)
For the $\vecv$-variables we impose again the initial 
datum $(\ref{ICv})$.

\begin{figure}
\vskip 3mm
\centerline{
\hskip 3mm
\includegraphics[width=0.48\columnwidth]{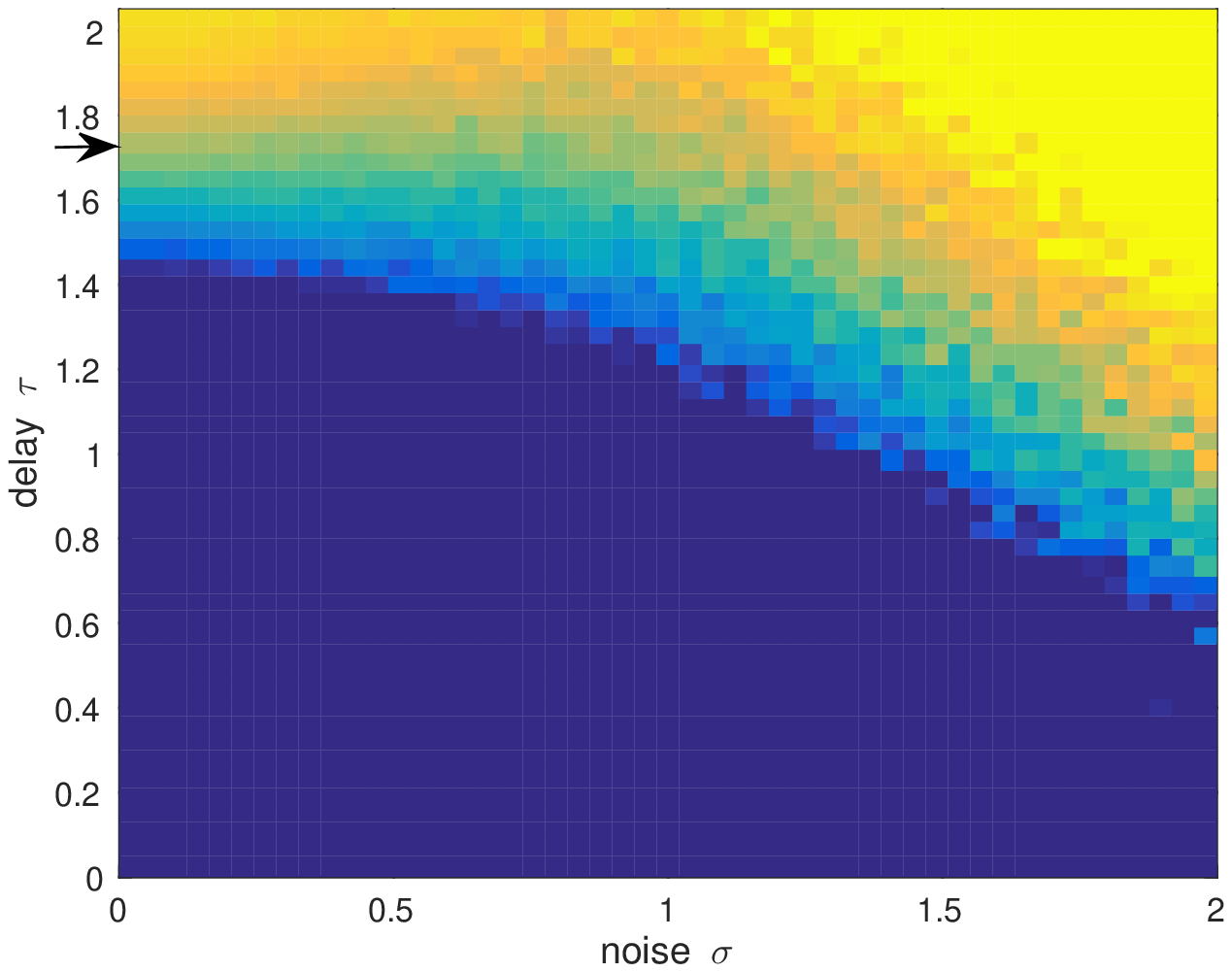}
\hskip 3mm
\includegraphics[width=0.48\columnwidth]{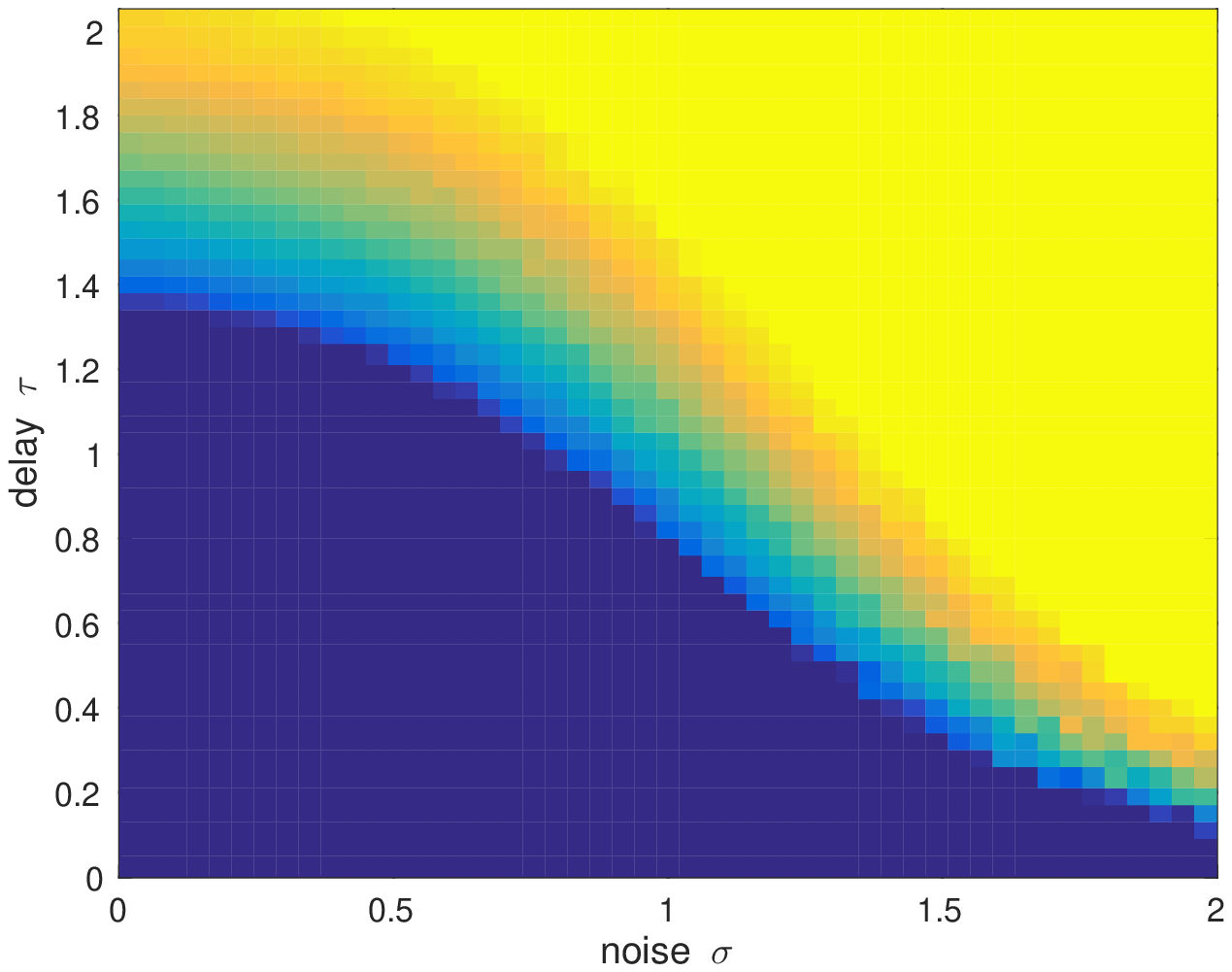}
}
\vskip -5cm
\leftline{\hskip 2mm (a) $\beta=0.1$, $N=2$ \hskip 3.65cm 
(b) $\beta=0.1$, $N=20$}
\vskip 4.5cm
\vskip 3mm
\centerline{
\hskip 3mm
\includegraphics[width=0.48\columnwidth]{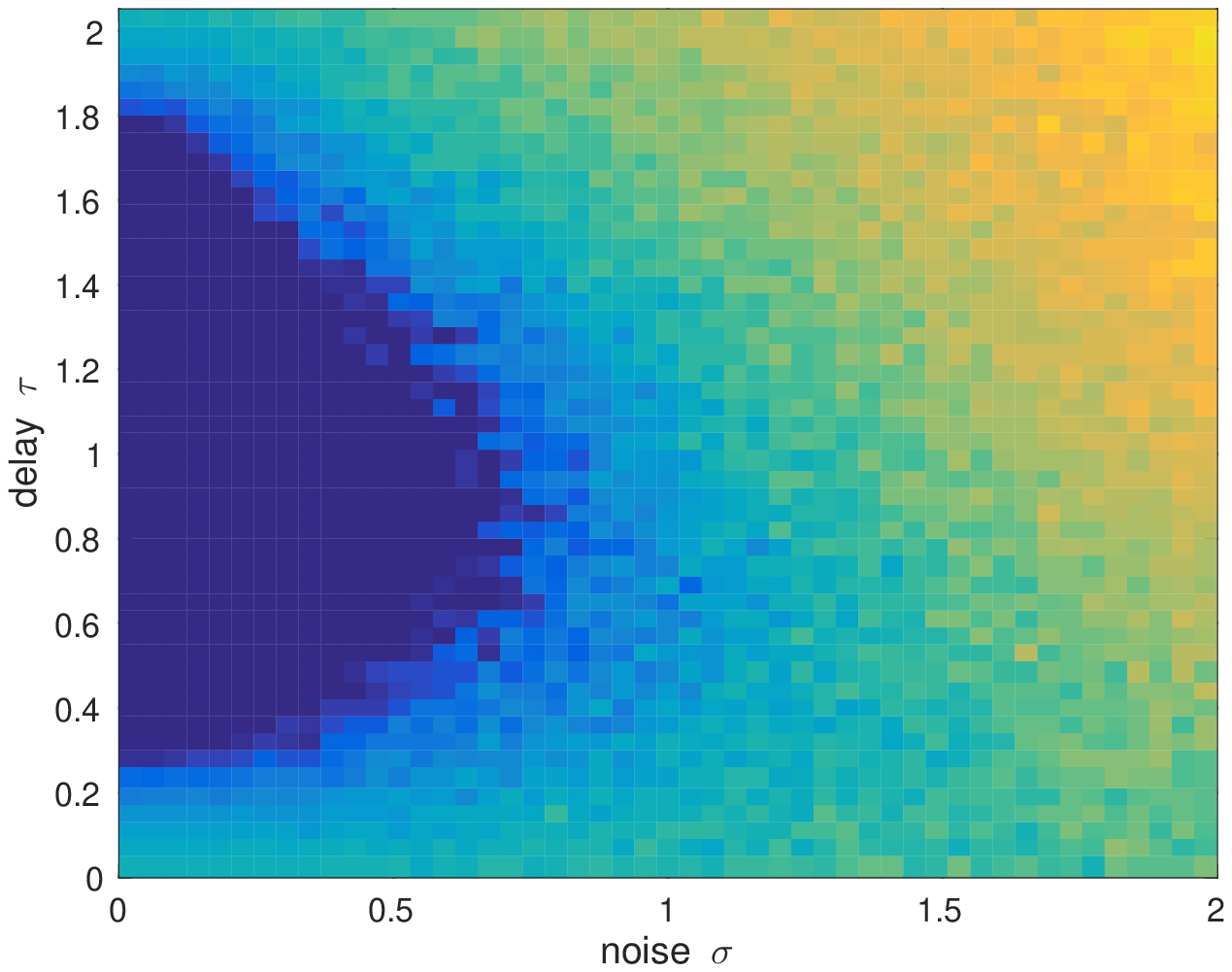}
\hskip 3mm
\includegraphics[width=0.48\columnwidth]{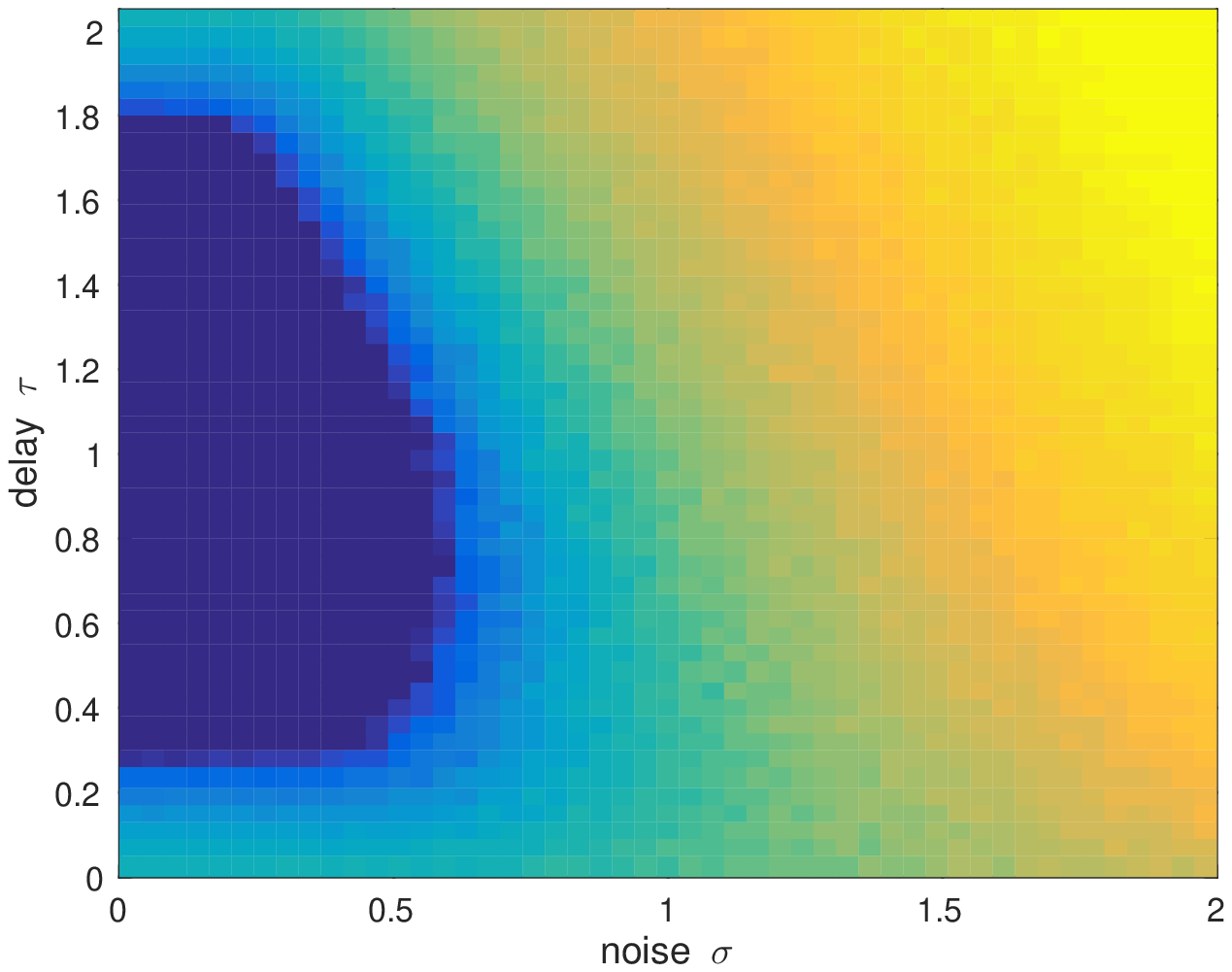}
}
\vskip -5cm
\leftline{\hskip 2mm (c) $\beta=1$, $N=2$ 
\hskip 4.05cm (d) $\beta=1$, $N=20$}
\vskip 4.5cm
\centerline{
\resizebox*{0.6\linewidth}{!}{\includegraphics{./colorbar_horz.eps}}
}
\caption{
{\it
Decadic logarithm of the indicator $I_{\sigma,\tau}$, given by
\eqref{indicator},
for Monte Carlo simulations of the system \eqref{CS01}--\eqref{CS02},
\eqref{commRate} with $Q=100$
paths on the time interval $[0,30]$ subject 
to the initial conditions \eqref{ICv} and \eqref{ICx},
with $\lambda=1$ and $(\sigma,\tau)\in[0,2]\times[0,2]$.
The dark blue regions (colour online) indicate numerical flocking.
The arrow in the top left panel indicates the 
point $(\sigma,\tau)=(0,1.75)$ that
corresponds to the parameter setting in Figure~$\ref{fig:IC}$.
We use:}
{\rm (a)} $\beta=0.1$, $N=2$; 
{\rm (b)} $\beta=0.1$, $N=20$;
{\rm (c)} $\beta=1$, $N=2$; {\it and} 
{\rm (d)} $\beta=1$, $N=20$.
\label{fig:CS-syst}}
\end{figure}

We perform Monte Carlo simulations of the system 
\eqref{CS01}--\eqref{CS02}, \eqref{commRate}
with $N\in\{2,20\}$ and $\beta=0.1$ (strong coupling) 
and $\beta=1$ (weak coupling).
As in Section~\ref{subsec:numCS0}, we fix $\lambda=1$ 
and divide the domain $[0,2]\times[0,2]$ for $(\sigma,\tau)$
into $50\times 50$ equidistant $(\sigma,\tau)$-pairs. For each 
pair of the parameter values we perform
a Monte Carlo simulation with $Q=100$ paths over the time interval
$[0,T]$ with $T=30$. We use the Euler-Maruyama method
for discretization of \eqref{CS01}--\eqref{CS02} with 
timestep $\Delta t = 10^{-3}$.
To classify the asymptotic behaviour of the solution, 
we again use the indicator
$(\ref{indicator})$ and say that \emph{numerical flocking} takes 
place when $I_{\sigma,\tau} < 10^{-2}$.
The background colour in Figure~\ref{fig:CS-syst} encodes 
the decadic logarithm of the indicator,
and the dark blue region indicates numerical flocking.

In the top left panel we indicate by an arrow the point 
$(\sigma,\tau)=(0,1.75)$
that corresponds to the parameter setting in Figure~\ref{fig:IC};
however, note that the initial conditions for $x_i$ in Figure~\ref{fig:IC}
differ from \eqref{ICx}.
We see that the indicated point lies close to the boundary of 
the dark blue region,
i.e., in the ``transition zone'' between numerical flocking and non-flocking.
We hypothesize that this is why we were able to observe the two qualitatively
different kinds of asymptotic behaviour in Figure~\ref{fig:IC}
even if the initial datum for the $\vecv$-variables is the same 
in both cases.
Again, a systematic study of this hypothesis is beyond the 
scope of this paper.

\begin{figure}
\vskip 3mm
\centerline{
\hskip 3mm
\includegraphics[width=0.48\columnwidth]{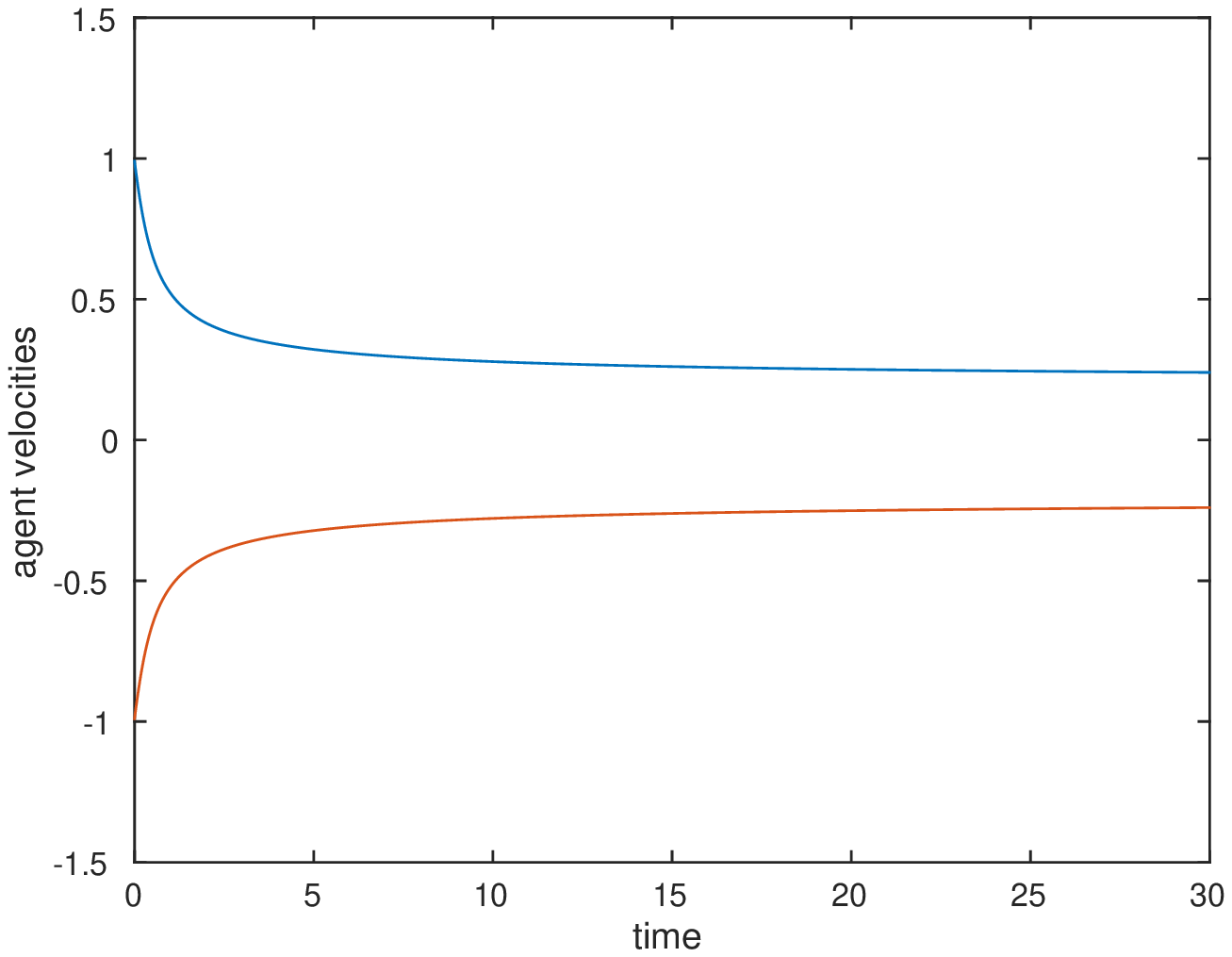}
\hskip 3mm
\includegraphics[width=0.48\columnwidth]{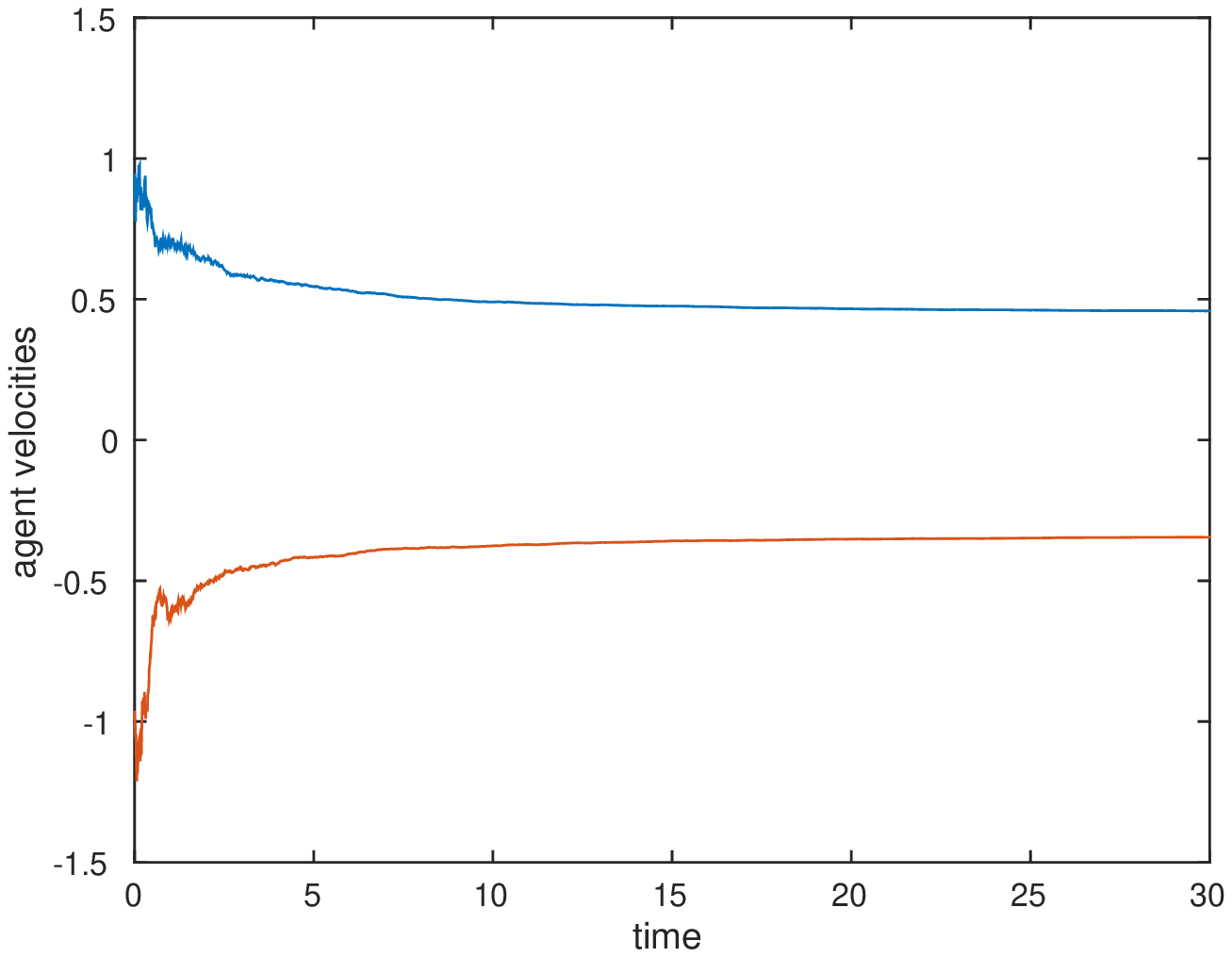}
}
\vskip -5cm
\leftline{\hskip 2mm (a) $\sigma=0$ \hskip 5.15cm (b) $\sigma=0.5$}
\vskip 3mm
\hskip 2.5cm $\tau=0$, $\sigma=0$ \hskip 4.6cm $\tau=0$, $\sigma=0.5$ 
\vskip 3.5cm
\centerline{
\hskip 3mm
\includegraphics[width=0.48\columnwidth]{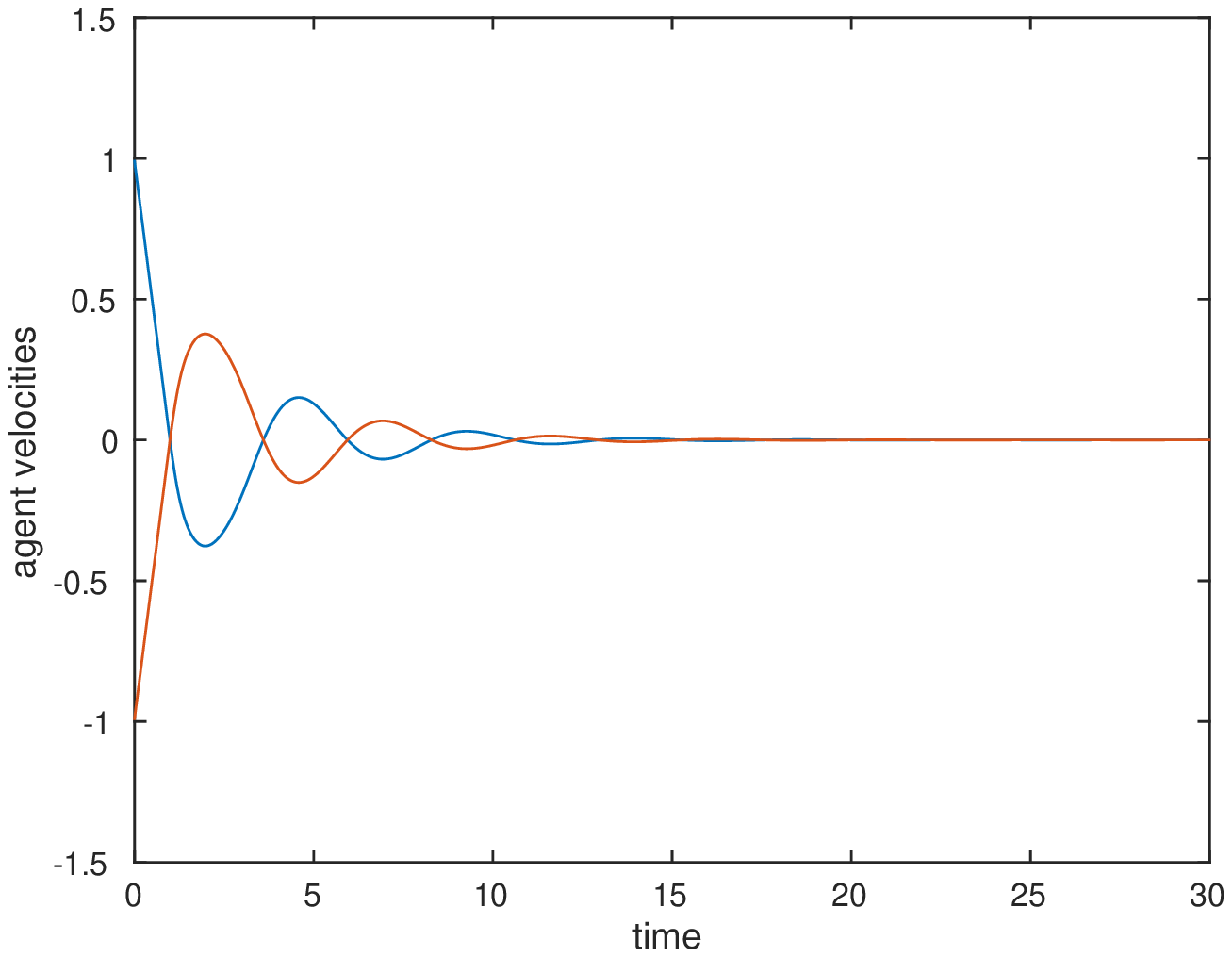}
\hskip 3mm
\includegraphics[width=0.48\columnwidth]{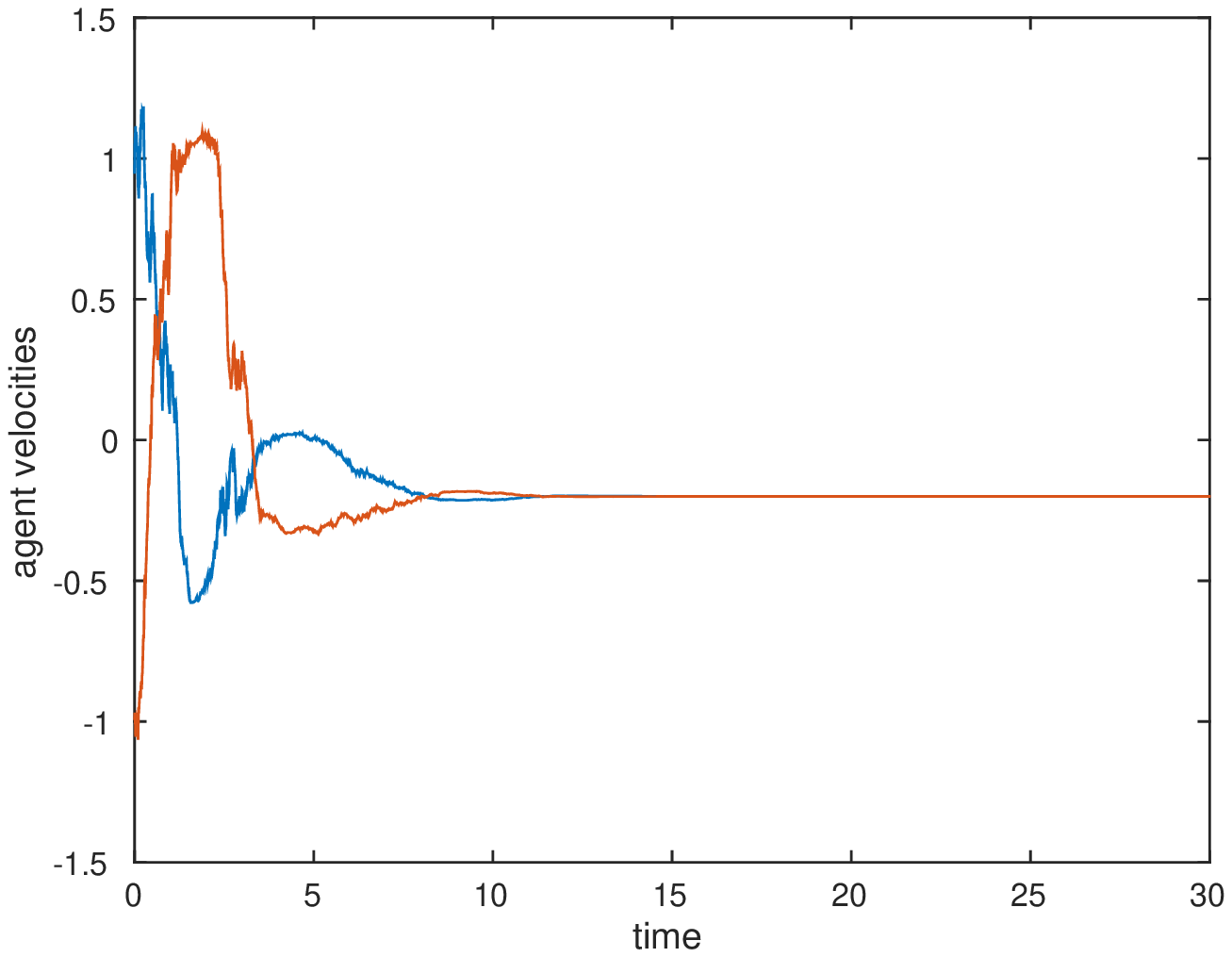}
}
\vskip -4.2cm
\hskip 2.5cm $\tau=1$, $\sigma=0$ \hskip 4.6cm $\tau=1$, $\sigma=0.5$ 
\vskip 3.5cm
\centerline{
\hskip 3mm
\includegraphics[width=0.48\columnwidth]{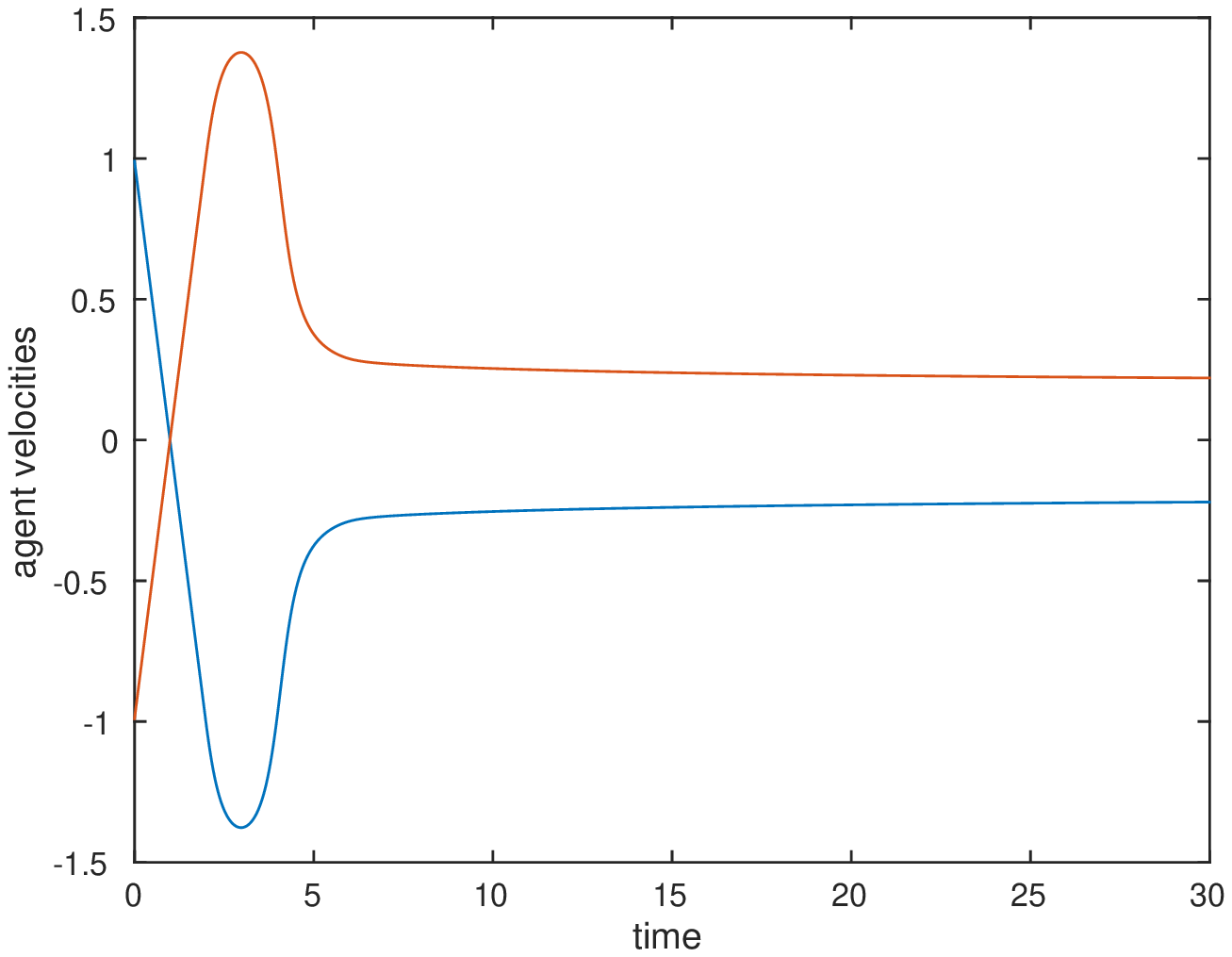}
\hskip 3mm
\includegraphics[width=0.48\columnwidth]{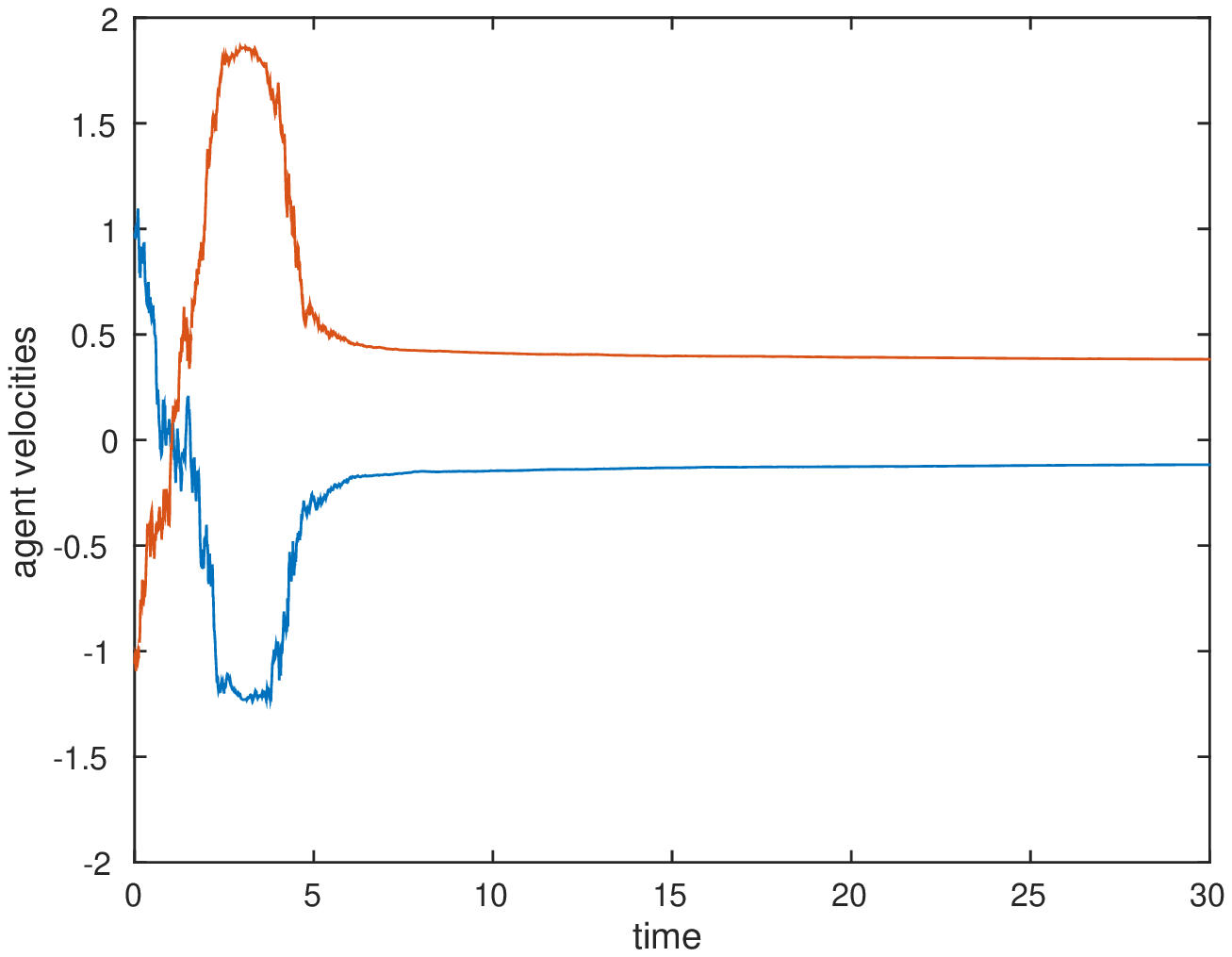}
}
\vskip -4.2cm
\hskip 2.5cm $\tau=2$, $\sigma=0$ \hskip 4.6cm $\tau=2$, $\sigma=0.5$ 
\vskip 3.5cm
\caption{
{\it Agents velocities $v_1(t), v_2(t)$ in sample solutions of the 
system \eqref{CS01}--\eqref{CS02}, \eqref{commRate}
with $N=2$, $\lambda=1$, $\beta=1$ (weak coupling),
on the time interval $[0,30]$ subject to the initial 
conditions \eqref{ICv} and \eqref{ICx}. We use:}
{\rm (a)} $\sigma=0$; {\it and} {\rm (b)} $\sigma=0.5$.
\label{fig:beta-tau-ex}}
\end{figure}

In Figure~\ref{fig:CS-syst} we observe that the region 
of numerical flocking is only weakly influenced by the number of agents $N$.
This is in agreement with the fact that the flocking 
condition \eqref{flockingCond} in Theorem \ref{thm:flocking}
does not depend on $N$.
The increased smoothness of the colour transition when $N=20$ 
is a consequence of the law of large numbers.
On the other hand, we can distinguish two distinct types of patterns,
one similar to Figure \ref{fig:noise-delay} for the strong coupling 
case $\beta=0.1$ (Figures~\ref{fig:CS-syst}(a) 
and~\ref{fig:CS-syst}(b)), and a semicircular 
pattern for the weak coupling case $\beta=1$ 
(Figures~\ref{fig:CS-syst}(c) and~\ref{fig:CS-syst}(d)).
In particular, the result for the weak coupling case is 
somewhat surprising -- it suggests 
that for low levels of noise ($\sigma\lesssim 0.6$), introduction 
of intermediate delays
($0.3\lesssim \tau \lesssim 1.8$) may facilitate flocking.
This is further supported by Figure \ref{fig:beta-tau-ex} where we plot
sample solutions of \eqref{CS01}--\eqref{CS02}, \eqref{commRate} 
for $N=2$, $\beta=1$, $\sigma=0$ (Figure~\ref{fig:beta-tau-ex}(a)), 
$\sigma=0.5$ (Figure~\ref{fig:beta-tau-ex}(b))
and three different values of the delay $\tau\in\{0,1,2\}$.
We observe that while for $\tau=0$ and $\tau=2$ the agents
do not show tendency to converge to a common velocity during the 
indicated time interval,
they exhibit numerical flocking for the intermediate value $\tau=1$.
We will call this observation \emph{time-delay induced flocking}.

Let us note that the results presented in Figure \ref{fig:beta-tau-ex} 
do not contradict our analytical results.
In particular, condition \eqref{flockingCond} gives
$\tau < \sqrt{2}/4 \doteq 0.35$ if $\sigma=0$
and $\tau < (-1/2+\sqrt{11/8})/4 \doteq 0.17$ if $\sigma=0.5$,
so it is only satisfied for the simulations in the panels
corresponding to $\tau=0$ in Figure~\ref{fig:beta-tau-ex}.
Therefore, statement of Lemma \ref{lem:Lyapunov} applies. 
The (expectation of)
the Lyapunov function \eqref{Lyapunov} decreases in time for 
these two simulations.

\begin{figure}
\vskip 3mm
\centerline{
\hskip 3mm
\includegraphics[width=0.48\columnwidth]{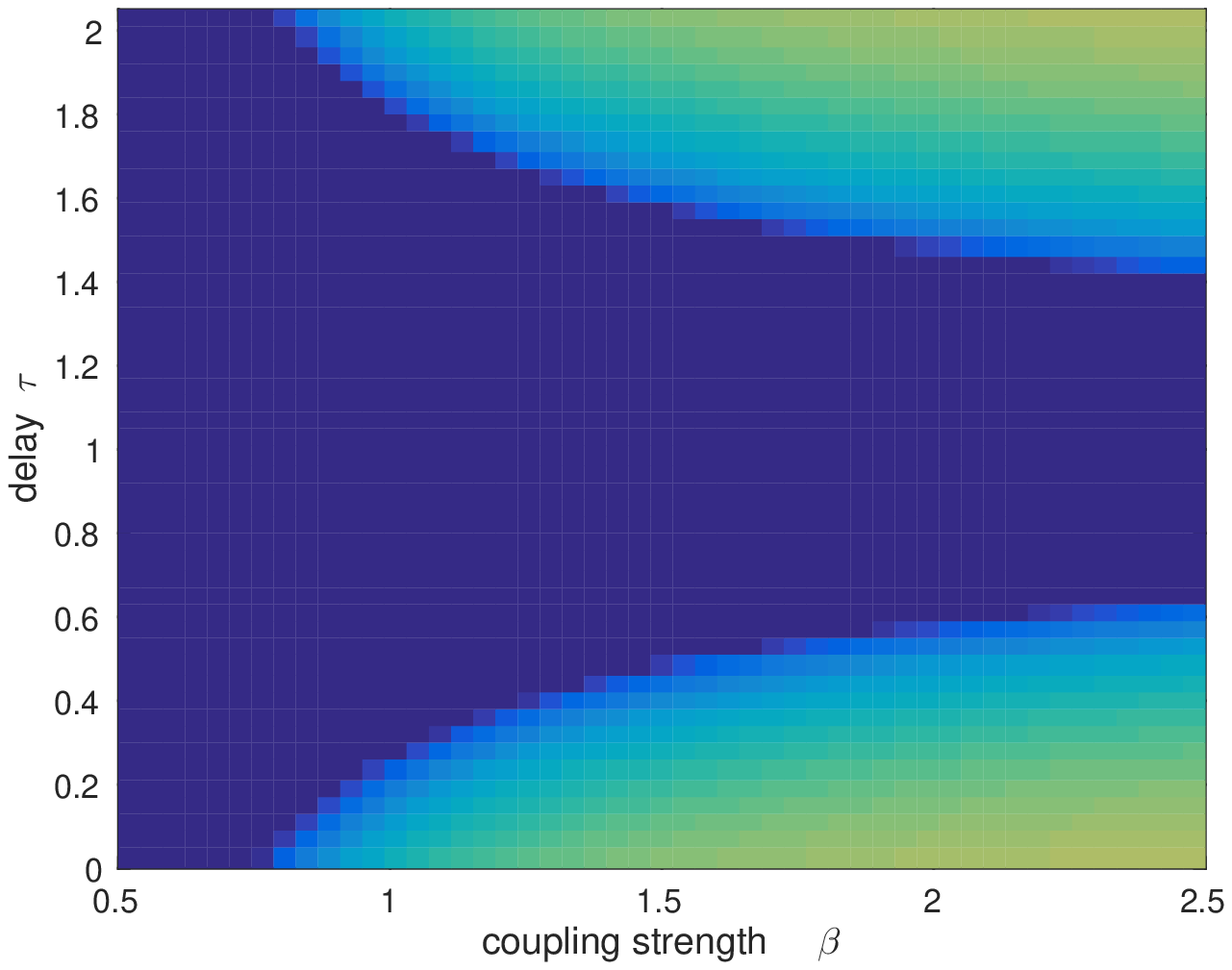}
\hskip 3mm
\includegraphics[width=0.48\columnwidth]{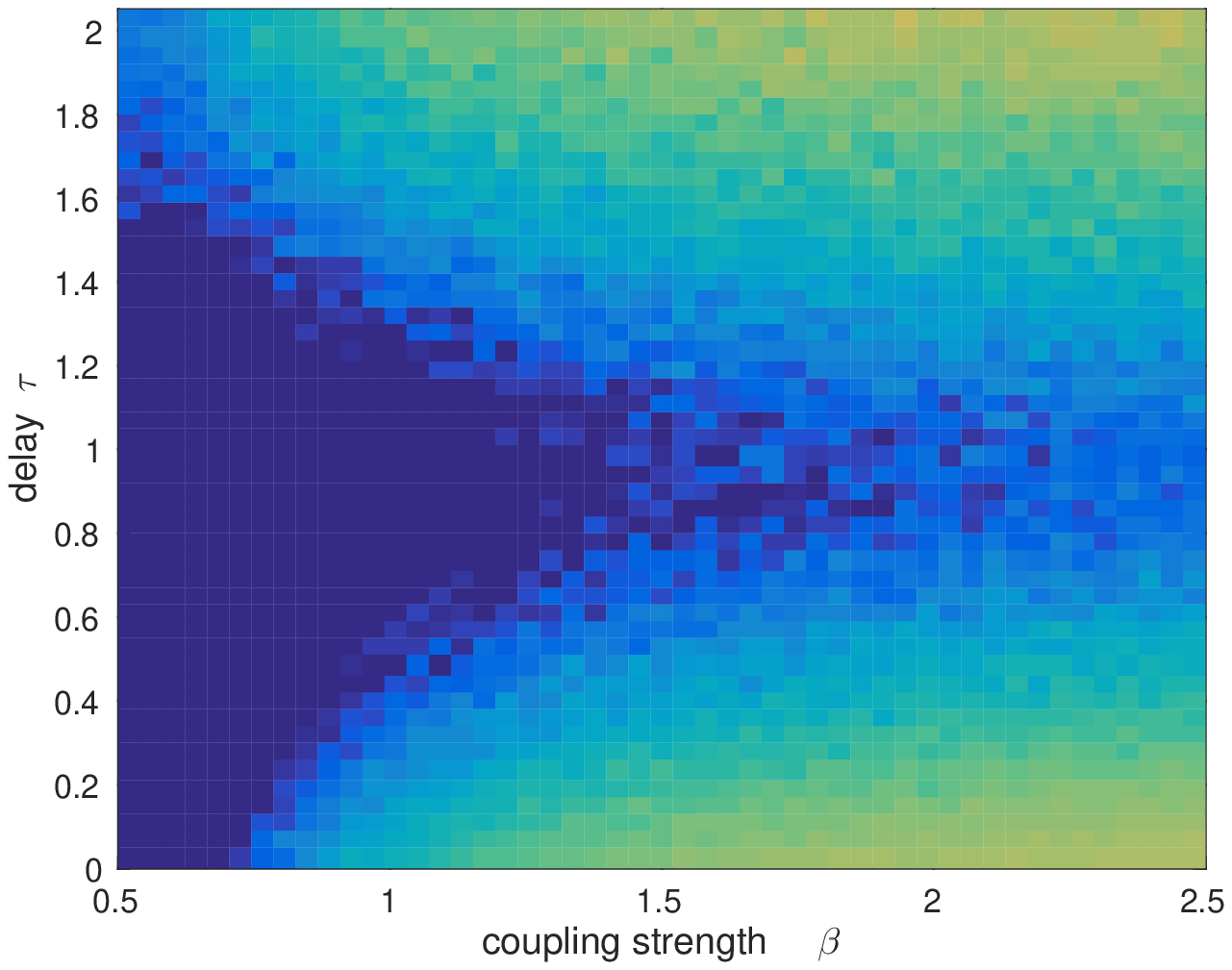}
}
\vskip -5cm
\leftline{\hskip 2mm (a) $\sigma=0$ \hskip 5.15cm (b) $\sigma=0.5$}
\vskip 4.5cm
\centerline{
\resizebox*{0.6\linewidth}{!}{\includegraphics{./colorbar_horz.eps}}
}
\caption{{Decadic logarithm of the indicator $I_{\beta,\tau}$
for simulations of the system \eqref{CS01}--\eqref{CS02}, 
\eqref{commRate} on the time interval $[0,30]$ with 
$\lambda=1$, $N=2$, $(\beta,\tau)\in[0.5,2.5]\times[0,2]$.
We use:} {\rm (a)} $\sigma=0$, $Q=1$; and
{\rm (b)} $\sigma=0.5$, $Q=100$.
{\it
The dark blue regions (colour online) indicate numerical flocking.}
\label{fig:beta-tau}}
\end{figure}

To gain a further understanding of the interesting phenomenon of 
time-delay induced flocking,
we run systematic simulations of the system 
\eqref{CS01}--\eqref{CS02}, \eqref{commRate}
with different values of $\beta\in [0.5,2.5]$, $\tau\in [0,2]$
and $\sigma\in\{0, 0.5\}$. We calculate the indicator $I_{\beta,\tau}$
as in \eqref{indicator} with $Q=1$ for $\sigma=0$ 
(there is no need to run more than one path for the case
without noise) and $Q=100$ Monte Carlo paths for $\sigma=0.5$.
The decadic logarithm of $I_{\beta,\tau}$ is plotted 
in Figure~\ref{fig:beta-tau}
and we again use the threshold $I_{\beta,\tau} < 10^{-2}$
to define numerical flocking (dark blue regions in Figure~\ref{fig:beta-tau}).
We observe that there exists
(for $\beta$ sufficiently large) a region of intermediate
values of $\tau$ where numerical flocking takes place,
while it does not for smaller or larger $\tau$ values.
Moreover, we see that noise has a disruptive influence on flocking
(the dark blue region is smaller in Figure~\ref{fig:beta-tau}(b) compared 
to Figure~\ref{fig:beta-tau}(a)).

\section{Discussion}

We have studied a generalization of the Cucker-Smale model
accounting for measurement errors through introduction
of multiplicative white noise, and for delays in information processing.
This has led to a system of stochastic delayed differential
equations~\eqref{CS01}--\eqref{CS02}.
In Section~\ref{sec:Flocking}, we have considered the communication rates
between agents as given stochastic processes, and derived a sufficient
condition for \emph{flocking}, which we define as asymptotic convergence 
of the agents' velocities towards a common value. The condition is given 
in terms of the critical delay that guarantees flocking as a function
of the noise level. Our analysis is based on a construction
of a suitable Lyapunov function for the system and a study of its decay.
As a byproduct of the analysis, we obtain a sufficient condition
for asymptotic convergence of delayed geometric Brownian motion.

The second part of the paper is devoted to systematic
numerical simulations. First, we perform Monte Carlo simulations
of delayed geometric Brownian motion and evaluate its asymptotic
behaviour based on a suitable ``numerical indicator''.
This led to the conclusion that the analytically derived sufficient 
condition for asymptotic convergence is qualitatively right - 
the convergence deteriorates with increasing noise level and delay.
However, quantitatively it is far from optimal.
Next, we simulate the Cucker-Smale type system with fixed communication rates
and again compare with the analytical result. As before, the comparison
shows that, while qualitatively correct, the analytical formula
produces too restrictive critical delays.
Finally, we simulate the full Cucker-Smale system with delays and
multiplicative noise. We use two regimes for the dependence of 
the communication rates
on the agents' distances: the strong coupling regime, which leads to
\emph{unconditional flocking} in the ``classical'' Cucker-Smale model,
and the weak coupling regime, where flocking may or may not take place.
In the strong coupling regime the numerical picture
is similar to the previous simulation with
fixed communication rates. On the other hand,
in the weak coupling regime we observe a somehow surprising behaviour
of the system - namely, that an introduction of intermediate time delay
may facilitate flocking. We call this phenomenon ``delay induced flocking''.

Our paper leaves several open questions.
First of all, our analytical flocking condition is too restrictive
compared to numerical results, so efforts should be made to improve it.
Moreover, the analysis applied to the case when the communication 
rates are given and satisfying a certain structural assumption.
This is in fact against the spirit of the original Cucker-Smale
model where the communication rates depend on the mutual distances
between agents. A possible extension of our analysis to this case remains
an open problem. The main difficulty is due to the fact that it is not
clear how to apply the classical bootstrapping argument that bounds
the velocity fluctuations in terms of fluctuations in positions and 
vice versa. For the numerical part, it would be desirable to apply 
some multilevel Monte Carlo
or importance sampling technique to obtain more accurate results.
Moreover, the influence of the initial condition on the asymptotic 
behaviour should be studied. Finally, the interesting phenomenon 
of ``delay induced flocking''
deserves a detailed study, both from the analytical and numerical
point of view.

\end{document}